\documentclass{tac}

\usepackage{amsmath}
\usepackage{amssymb}
\usepackage[pdf]{xy}
\usepackage{stmaryrd}

\usepackage{hyperref}
\hypersetup{colorlinks,%
citecolor=black,%
filecolor=black,%
linkcolor=black,%
urlcolor=black,%
pdftex}

\input diagxy
\xyoption {all}
\xyoption{2cell}
\xyoption{v2}

\author {David Michael Roberts}

\thanks{An Australian Postgraduate Award provided financial support during
part of the time the material for this paper was written.}

\address{School of Mathematical Sciences,\\ University of Adelaide\\
Adelaide, SA 5005\\Australia
}

\title {Internal categories, anafunctors and localisations}

\copyrightyear{2011}

\eaddress{david.roberts@adelaide.edu.au}

\dedication{In memory of Luanne Palmer (1965-2011)}

\newcounter{dr}

\newtheorem{theorem}{Theorem}

\newtheorem{proposition}{Proposition}
\newtheorem{lemma}{Lemma}
\newtheorem{corollary}{Corollary}

\newtheoremrm{remark}{Remark}
\newtheoremrm{definition}{Definition}
\newtheoremrm{example}{Example}

\let\endofproof\endproof

\mathrmdef{Hom}
\mathrmdef{Arr}
\mathrmdef{Mor}
\mathrmdef{Obj}
\mathrmdef{ana}

\mathbfdef{Set}
\mathbfdef{Ab}
\mathbfdef{Top}
\mathbfdef{Ring}
\mathbfdef{Diff}
\mathbfdef{Aff}
\mathbfdef{Grp}
\mathbfdef{Cat}
\mathbfdef{Sch}
\mathbfdef{Gpd}
\mathbfdef{Lie}
\mathbfdef{Fields}
\mathbfdef{AlgSp}

\mathopdef{colim}
\mathopdef{im}
\mathopdef{coker}
\mathopdef{op}
\mathopdef{id}
\mathopdef{pr}
\mathopdef{Aut}
\mathopdef{codisc}
\mathopdef{disc}

\pdfoutput=1 

\usepackage{microtype}

\newcommand{\st}{\rightrightarrows} 
\newcommand{\into}{\hookrightarrow}    

\newcommand{\can}{\underline{c}}

\newcommand{\gento}{-\!\!\!\mapsto}
\newcommand{\equals}{\hspace{10pt}=\hspace{10pt}}   

\begin{document}

   \maketitle

	\begin{abstract}
	
	In this article we review the theory of anafunctors introduced by Makkai
	and Bartels, and show that given a subcanonical site $S$, one can form a 
	bicategorical localisation of various 2-categories of internal categories or
	groupoids at weak equivalences using anafunctors as 1-arrows. This unifies 
	a number of proofs throughout the literature, using the fewest assumptions 
	possible on $S$.

	\end{abstract}


\section{Introduction}\label{section_1}

It is a well-known classical result of category theory that a functor
is an equivalence (that is, in the 2-category of categories) if and
only if it is fully faithful and essentially surjective.
This fact is equivalent to the axiom of choice. It is therefore
\emph{not} true if one is working with categories internal to a category
$S$ which doesn't satisfy the (external) axiom of choice. This is may fail
even in a category very much like the category of sets, such as a 
well-pointed boolean topos, or even the category of sets in constructive
foundations.  As internal
categories are the objects of a 2-category $\Cat(S)$ we can talk about
internal equivalences, and even fully faithful functors. In the case
$S$ has a singleton
pretopology $J$ (i.e.~covering families 
consist of single maps) we can define an analogue of essentially
surjective functors. Internal functors which are fully faithful
and essentially surjective are called \emph{weak equivalences} in the
literature, going back to \cite{Bunge-Pare_79}. We shall call them $J$-equivalences for
clarity. We can recover the classical result mentioned above if we localise the
2-category $\Cat(S)$ at the class $W_J$ of $J$-equivalences.

We are not just interested in localising $\Cat(S)$, but various
full sub-2-categories $C \into \Cat(S)$ which arise in the
study of presentable stacks, for example algebraic, topological, differentiable,
etc.~stacks. As such it is necessary to ask for a compatibility
condition between the pretopology on $S$ and the sub-2-category we are
interested in. We call this condition existence of \emph{base change} for
covers of the pretoplogy, and demand that for any cover $p\colon U\to X_0$
(in $S$) of the object of objects of $X\in C$, there is a fully faithful
functor in $C$ with object component $p$.

\begin{theorem}
Let $S$ be a category with singleton pretopology $J$ and let $C$ be a
full sub-2-category of $\Cat(S)$ which admits base change along arrows in
$J$. Then $C$ admits a calculus of fractions for the $J$-equivalences.
\end{theorem}

Pronk gives us the appropriate notion of a calculus of fractions for a
2-category in \cite{Pronk_96} as a generalisation of the usual construction
for categories \cite{Gabriel-Zisman}. In her construction, 1-arrows are spans and 2-arrows are
equivalence classes of bicategorical spans of spans. This construction,
while canonical, can be a little unwieldy so we look for a simpler construction
of the localisation.

We find this in the notion of \emph{anafunctor}, introduced by Makkai
for plain small categories \cite{Makkai} (Kelly described them briefly
in \cite{Kelly_64} but did not develop the concept further). In his setting an anafunctor
is a span of functors such that the left (or source) leg is a surjective-on-objects,
fully faithful functor.\footnote{Anafunctors were so named by Makkai, on the 
suggestion of Pavlovic, after profunctors,
in analogy with the pair of terms anaphase/prophase from biology. For more on the 
relationship between anafunctors and profunctors, see below.}
For a general category $S$ with a \emph{subcanonical}
singleton pretopology $J$ \cite{Bartels}, the analogon is a span with left leg a 
fully faithful functor with object component a cover. Composition of 
anafunctors is given by composition of spans in the usual way, and there 
are 2-arrows between anafunctors (a certain sort of span of spans) that 
give us a bicategory $\Cat_\ana(S,J)$ with objects internal categories
and 1-arrows anafunctors. We can also define the full sub-bicategory
$C_\ana(J) \into \Cat_\ana(S,J)$ analogous to $C$, and there is a strict
inclusion 2-functor $C \into C_\ana(J)$. This gives us our second main
theorem.

\begin{theorem}
Let $S$ be a category with subcanonical singleton pretopology $J$ and let
$C$ be a full sub-2-category of $\Cat(S)$ which admits base change
along arrows in $J$, Then $C \into C_\ana(J)$ is a localisation
of $C$ at the class of $J$-equivalences.
\end{theorem}

So far we haven't mentioned the issue of size, which usually is important
when constructing localisations. If the site $(S,J)$ is locally small,
then $C$ is locally small, in the sense that the hom-categories are
small. This also implies that $C_\ana(J)$ and hence any 
$C[W_J^{-1}]$ has \emph{locally} small hom-categories i.e.~has only
a set of 2-arrows between any pair of 1-arrows. To prove that the
localisation is locally essentially small (that is, hom-categories are equivalent
to small categories), we need to assume a size restriction axiom on the 
pretopology $J$, called WISC (Weakly Initial Sets of Covers).

WISC can be seen as an extremely weak choice principle, weaker than the existence
of enough projectives, and states that for every object $A$ of $S$, there is
a set of $J$-covers of $A$ which is cofinal in all $J$-covers of $A$. It is 
automatically satisfied if the pretopology is specified as an assignment of
a \emph{set} of covers to each object.

\begin{theorem}
Let $S$ be a category with subcanonical singleton pretopology $J$ satisfying
WISC, and let $C$ be a full sub-2-category of $\Cat(S)$ which admits
base change along arrows in $J$. Then any localisation of $C$ at the 
class of $J$-equivalences is locally essentially small.
\end{theorem}

Since a singleton pretopology can be conveniently defined as a certain wide 
subcategory, this is not a vacuous statement for large sites, such as $\Top$
or $\Grp(E)$ (group objects in a topos $E$). In fact WISC is independent
of the Zermelo-Fraenkel axioms (without Choice) \cite{vdBerg_12,Roberts_13}.
It is thus possible to have the theorem fail for the topos $S = \Set_{\neg AC}$ 
with surjections as covers.

Since there have been many very closely related approaches to localisation
of 2-categories of internal categories and groupoids, we give a brief sketch 
in the following section. Sections 3 to 6 of this article then give necessary background and notation
on sites, internal categories, anafunctors and bicategories of fractions respectively. 
Section 7 contains our main results, while section 8 shows examples from
the literature that are covered by the theorems from section 7. A short appendix
detailing superextensive sites is included, as this material does not appear to
be well-known (they were discussed in the recent \cite{Shulman_12}, Example 11.12).

This article started out based on the first chapter of the author's PhD thesis,
which only dealt with groupoids in the site of topological spaces
and open covers.
Many thanks are due to Michael Murray, Mathai Varghese and Jim Stasheff,
supervisors to the author. The patrons of the $n$-Category Caf\'e and $n$Lab,
especially Mike Shulman and Toby Bartels, provided helpful input and feedback.
Steve Lack suggested a number of improvements, and the referee asked for
a complete rewrite of this article, which
has greatly improved the theorems, proofs, and hopefully also the exposition. Any delays
in publication are due entirely to the author.

\section{Anafunctors in context}

The theme of giving 2-categories of internal categories or groupoids more
equivalences has been approached in several different ways over the decades. 
We sketch a few of them, without necessarily finding the original references, to
give an idea of how widely the results of this paper apply. We give some 
more detailed examples of this applicability in section 8.
 
Perhaps the oldest related construction is the distributors of B\'enabou, also
known as modules or profunctors \cite{Benabou_73} (see \cite{Elephant} for 
a detailed treatment of internal profunctors, as the original article is difficult to source). 
B\'enabou pointed out \cite{Benabou_email}, after a preprint of this article was 
released, that in the case of the category $\Set$ (and more generally in a finitely 
complete site with reflexive coequalisers that are stable under pullback, 
see \cite{MMV2012}), the bicategory
of small (resp. internal) categories with representable profunctors as 1-arrows is 
equivalent to the bicategory of small categories with anafunctors as 1-arrows.
In fact this was discussed by Baez and Makkai \cite{Baez-Makkai_emails}, 
where the latter pointed out that representable profunctors correspond to
\emph{saturated} anafunctors in his setting. The author's preference for anafunctors
lies in the fact they can be defined with weaker assumptions on the site $(S,J)$,
and in fact in the sequel \cite{Roberts2}, do not require the 2-category to have objects
which are internal categories. In a sense this is analogous to \cite{Street_80},
where the formal bicategorical approach to profunctors between objects of
a bicategory is given, albeit still requiring more colimits to exist than anafunctors do.

B\'enabou has pointed out in private communication that he has 
an unpublished distributor-like construction that does not rely on existence of 
reflexive coequalisers; the author has not seen any details of this and is 
curious to see how it compares to anafunctors.

Related to this is the original work of Bunge and Par\'e \cite{Bunge-Pare_79}, 
where they consider functors between indexed categories associated to 
internal categories, that is, the \emph{externalisation} of an internal category and stack
completions thereof. This was
one motivation for considering weak equivalences in the first place, in that a pair
of internal categories have equivalent stack completed externalisations if and only if they are 
connected by a span of internal functors which are weak equivalences.

Another approach is constructing bicategories of fractions \`a la Pronk \cite{Pronk_96}.
This has been followed by a number of authors, usually followed up by an explicit
construction of a localisation simplifying the canonical one. Our work here sits at the
more general end of this spectrum, as others have tailored their constructions to
take advantage of the structure of the site they are interested in. For example, 
\emph{butterflies} (originally called papillons) have been used for the category of 
groups \cite{Noohi_05b, Aldrovandi-Noohi_09,Aldrovandi-Noohi_10},
abelian categories \cite{Breckes_09} and semiabelian categories \cite{AMMV_10,MMV2012}. 
These are similar to the meromorphisms of \cite{Pradines_89}, introduced in 
the context of the site of smooth manifolds; though these only use a 1-categorical 
approach to localisation. 

Vitale \cite{Vitale_10}, after first showing that the 
2-category of groupoids in a regular category has a bicategory
of fractions, then shows that for protomodular regular categories one can
generalise the pullback congruences of B\'enabou in \cite{Benabou_89} to
discuss bicategorical localisation. This approach can be applied to internal
categories, as long as one restricts to invertible 2-arrows. Similarly, \cite{MMV2012}
give a construction of what they call \emph{fractors} between internal 
groupoids in a Mal'tsev category, and show that in an efficiently regular category 
(e.g.~a Barr-exact category) fractors are 1-arrows in a localisation of the 2-category of 
internal groupoids. The proof also works for internal categories if one considers
only invertible 2-arrows.

Other authors, in dealing with internal groupoids, have adopted the 
approach pioneered by Hilsum and Skandalis
\cite{Hilsum-Skandalis_87}, which has gone by various names including 
Hilsum-Skandalis morphisms, Morita morphisms, bimodules, bibundles, right principal bibundles
and so on. All of these are very closely related to saturated anafunctors, but
in fact no published definition of a saturated anafunctor in a site other than $\Set$
(\cite{Makkai}) has appeared, except in the guise of internal profunctors 
(e.g. \cite{Elephant}, section B2.7). Note also that this approach has only been applied to internal
groupoids. The review \cite{Lerman_10} covers the case of Lie groupoids, and in
particular orbifolds, while \cite{Mrcun_01} treats bimodules between groupoids in the 
category of affine schemes, but from the point of view of Hopf algebroids.

The link between localisation at weak equivalences and presentable stacks is considered
in (of course) \cite{Pronk_96}, as well as more recently in \cite{Carchedi_12}, \cite{Schappi_12},
in the cases of topological and algebraic stacks respectively, and for example
\cite{Tu-Xu-LaurentGengoux_04} in the case of differentiable stacks.

A third approach is by considering a model category structure on the 1-category
of internal categories. This is considered in \cite{Joyal-Tierney_91} for categories
in a topos, and in \cite{Everaert_et_al_05} for categories in a  
finitely complete subcanonical site $(S,J)$.
In the latter case the authors show when it is possible to construct a Quillen model category 
structure on $\Cat(S)$ where the weak equivalences are the weak equivalences
from this paper. Sufficient conditions on $S$ include being a topos with nno, being
locally finitely presentable or being finitely complete regular Mal'tsev -- and additionally
having enough $J$-projective objects. If one is willing to consider other model-category-like
structures, then these assumptions can be dropped. The proof from \cite{Everaert_et_al_05}
can be adapted to show that for a finitely complete site $(S,J)$, the category of groupoids
with source and target maps restricted to be $J$-covers has the structure of a category of fibrant objects,
with the same weak equivalences. We note that \cite{Colman-Costoya_09} gives a Quillen
model structure for the category of orbifolds, which are there defined to be proper topological
groupoids with discrete hom-spaces.

In a similar vein, one could consider a localisation using \emph{hammock} localisation
\cite{Dwyer-Kan_80a
} of a category of internal categories,
which puts one squarely in the realm of $(\infty,1)$-categories. Alternatively, one could
work with the $(\infty,1)$-category arising from a 2-category of internal categories, functors
and natural \emph{isomorphisms} and consider a localisation of this as given in, say 
\cite{Lurie_HTT}. However, to deal with general 2-categories of internal categories in this way,
one needs to pass to $(\infty,2)$-categories to handle the non-invertible 2-arrows. The 
theory here is not so well-developed, however, and one could see the results 
of the current paper as giving 
toy examples with which one could work. This is one motivation for making sure the
results shown in this paper apply to not just 2-categories of groupoids. Another is extending
the theory of presentable stacks from stacks of groupoids to stacks of categories \cite{Roberts1}.

\section{Sites}\label{sites_categories}

The idea of \emph{surjectivity} is a necessary ingredient when talking about
equivalences of categories---in the guise of just essential surjectivity---but it doesn't
generalise in a straightforward way from the category $\Set$. The necessary properties
of the class of surjective maps are encoded in the definition of a Grothendieck pretopology,
in particular a singleton pretopology. This section gathers definitions and
notations for later use.

\begin{definition}
A \emph{Grothendieck pretopology} (or simply \emph{pretopology}) on a category $S
$ is a collection $J$ of families
\[
       \{ (U_i \to A)_{i\in I} \}_{A\in \Obj(S)}
\]
of morphisms for each object $A \in S$ satisfying the following properties
\begin{enumerate}
       \item
               $(A' \stackrel{\sim}{\to} A)$ is in $J$ for every isomorphism $A'\simeq A$.
       \item
               Given a map $B \to A$, for every $(U_i \to A)_{i\in I}$ in $J$ the pullbacks
$B \times_A A_i$ exist and $(B \times_A A_i \to B)_{i\in I}$ is in $J$.

       \item
               For every $(U_i \to A)_{i\in I}$ in $J$ and for a collection $(V_k^i \to
U_i)_{k\in K_i}$ from $J$ for each $i \in I$, the family of composites
               \[
                       (V_k^i \to A)_{k\in K_i,i\in I}
               \]
               are in $J$.

\end{enumerate}
Families in $J$ are called \emph{covering families}. We call a category $S$ equipped with a
pretopology $J$ a \emph{site}, denoted $(S,J)$ (note that often one sees
a site defined as a category equipped with a Grothendieck \emph{topology}).

The pretopology $J$ is called a \emph{singleton} pretopology if
every covering family consists of a single arrow $(U \to A)$. In this case a covering
family is called a \emph{cover} and we call $(S,J)$ a \emph{unary} site.
\end{definition}

Very often, one sees the definition of a pretopology as being an assignment of a \emph{set}
covering families to each object. We do not require this, as  
one can define a singleton pretopology as a subcategory with certain properties, 
and there is not necessarily then
a set of covers for each object. One example is the category of groups with surjective 
homomorphisms as covers. This distinction will be important later.

One thing we will require is that sites come with \emph{specified} pullbacks of 
covering families. If one does not mind applying the axiom of choice (resp.~axiom of choice
for classes) then any small site (resp.~large site) can be so equipped. But often sites that 
arise in practice have more or less canonical choices for pullbacks, such as the category
of ZF-sets.

\begin{example}
The prototypical example is the pretopology $\mathcal{O}$ on $\Top$, where a covering family is an open cover.
The class of numerable open covers (i.e.~those that admit a subordinate
partition of unity \cite{Dold_63}) also forms a pretopology on $\Top$. Much of traditional bundle theory
is carried out using this site; for example the Milnor classifying space classifies
bundles which are locally trivial over numerable covers.
\end{example}

\begin{definition}\label{defn:effective_cov_fam}
A covering family $(U_i \to A)_{i\in I} $ is called \emph{effective} if $A$ is the colimit
of the following diagram: the objects are the $U_i$ and the pullbacks $U_i \times_A
U_j$, and the arrows are the projections
\[
       U_i \leftarrow U_i \times_A U_j \to U_j.
\]
If the covering family consists of a single arrow $(U \to A)$, this is the same as saying
$U \to A$ is a regular epimorphism.
\end{definition}

\begin{definition}
A site is called \emph{subcanonical} if every covering family is effective.
\end{definition}

\begin{example}
On $\Top$, the usual pretopology $\mathcal{O}$ of opens,  the pretopology of numerable covers
and that of open surjections are subcanonical.
\end{example}

\begin{example}
In a regular category, the class of regular epimorphisms forms a subcanonical singleton
pretopology. 
\end{example}

In fact we can make the following definition.

\begin{definition}
For a category $S$, the largest class of regular epimorphisms of which all pullbacks exist,
and which is stable under pullback, is called the \emph{canonical singleton pretopology}
and denoted $\can$.
\end{definition}

This is a to be contrasted to the canonical \emph{topology} on a category, which consists
of covering sieves rather than covers. The canonical singleton pretopology is the largest
subcanonical singleton pretopology on a category.

\begin{definition}\label{defn:saturation}
Let $(S,J)$ be a site. An arrow $P \to A$ in $S$ is called a \emph{$J$-epimorphism}
 if there is a covering family $(U_i \to A)_{i\in I}$ and a lift
\[
       \xymatrix{
                & P \ar[d] \\
               U_i \ar@{-->}[ur] \ar[r] & A
       }
\]
for every $i \in I$. A $J$-epimorphism is called \emph{universal} if its pullback 
along an arbitrary map exists. We denote the singleton pretopology of universal $J$-epimorphisms 
by  $J_{un}$.
\end{definition}

This definition of $J$-epimorphism is equivalent to the definition in III.7.5 in 
\cite{MacLane-Moerdijk}. The
dotted maps in the above definition are called local sections, after the case of 
the usual open cover pretopology on $\Top$. If $J$ is a singleton pretopology, 
it is clear that $J \subset J_{un}$.

\begin{example}
The universal $\mathcal{O}$-epimorphisms for the pretopology $\mathcal{O}$ of 
open covers on $\Diff$ form $Subm$, the pretopology of surjective submersions.
\end{example}

\begin{example}\label{eg:split_epis}
In a finitely complete category the universal $triv$-epimorphisms are the split
epimorphisms, where $triv$ is the \emph{trivial pretopology} where all covering
families consist of a single isomorphism. In $\Set$ with the axiom of
choice there are all the epimorphisms.
\end{example}

Note that for a finitely complete site $(S,J)$, $J_{un}$ contains $triv_{un}$,
hence all the split epimirphisms.

Although we will not assume that all sites we consider are finitely complete,
results similar to ours have, and so in that case we can say a little more,
given stronger properties on the pretopology. 

\begin{definition}
A singleton pretopology $J$ is called \emph{saturated} if whenever the composite
$A \stackrel{h}{\to} B \stackrel{g}{\to} C$ is in $J$, then $g\in J$.
\end{definition}

The concept of a saturated pretopology was introduced by B\'enabou under the name 
\emph{calibration} \cite{Benabou_75a}. It follows from the definition that a saturated singleton 
pretopology contains the split epimorphisms (take $h$ to be a section of the epimorphism $g$).

\begin{example}
The canonical singleton pretopology $\can$ in a regular category (e.g.~a topos) is saturated.
\end{example}

\begin{example}
Given a pretopology $J$ on a finitely complete category, $J_{un}$ is saturated.
\end{example}

Sometimes a pretopology $J$ contains a smaller pretopology that still has enough
covers to compute the same $J$-epimorphisms.

\begin{definition}
If $J$ and $K$ are two singleton pretopologies with $J \subset K$, such that
$K \subset J_{un}$, then $J$ is said to be \emph{cofinal} in $K$.
\end{definition}

Clearly $J$ is cofinal in $J_{un}$ for any singleton pretopology $J$.

\begin{lemma}
If $J$ is cofinal in $K$, then $J_{un} = K_{un}$.
\end{lemma}

We have the following lemma, which is essentially proved in \cite{Elephant}, C2.1.6.

\begin{lemma}\label{subcanonical_goes_up_cofinal}
If a pretopology $J$ is subcanonical, then so any pretopology in which it is cofinal.
In particular, $J$ subcanonical implies $J_{un}$ subcanonical.
\end{lemma}

As mentioned earlier, one may be given a singleton pretopology such that each 
object has more than a set's worth of covers. If such a pretopology contains a
cofinal pretopology with set-many covers for each object, then we can pass to 
the smaller pretopology and recover the same results (in a way that
will be made precise later). In fact, we can get away 
with something weaker: one could ask only that the category of all covers of an
object (see definition \ref{cover_slice} below) has a set of weakly initial
objects, and such set may not form a pretopology.
This is the content of the axiom WISC below. We first give some more
precise definitions.

\begin{definition}
A category $C$ has a \emph{weakly initial set} $\mathcal{I}$ of objects if for every
object $A$ of $C$ there is an arrow $O\to A$ from some object $O\in \mathcal{I}$.
\end{definition}

For example the large category $\Fields$ of fields has a weakly initial set, consisting of 
the prime fields $\{\mathbb{Q},\mathbb{F}_p|p\textrm{ prime}\}$. 
To contrast, the category of sets with surjections for arrows doesn't have a
weakly initial set of objects. Every small category has a weakly initial set, namely its set of objects.

We pause only to remark that the statement of the adjoint functor theorem can be 
expressed in terms of weakly initial sets.

\begin{definition}\label{cover_slice}
Let $(S,J)$ be a site. For any object $A$, the \emph{category of covers of $A$},
denoted $J/A$ has as objects the covering families $(U_i \to A)_{i\in I}$ and as
morphisms $(U_i \to A)_{i\in I} \to (V_j \to A)_{j\in J}$ tuples consisting of a function
$r\colon I\to J$ and arrows $U_i \to V_{r(i)}$ in $S/A$.
\end{definition}
When $J$ is a singleton pretopology this is simply a full subcategory of $S/A$. 
We now define the axiom WISC (Weakly Initial Set of Covers), due independently to Mike
Shulman and Thomas Streicher.

\begin{definition}
A site $(S,J)$ is said to \emph{satisfy WISC} if for every object $A$ of $S$, the
category $J/A$ has a weakly initial set of objects.
\end{definition}

A site satisfying WISC is in some sense constrained by a small amount of data for each object.
Any small site satisfies WISC, for example, the usual site of finite-dimensional smooth
manifolds and open covers. Any pretopology $J$ containing a cofinal pretopology $K$
such that $K/A$ is small for every object $A$ satisfies WISC. 

\begin{example}
Any regular category (for example a topos) with enough projectives, equipped with the canonical
singleton pretopology, satisfies WISC. In the case of $\Set$ `enough projectives' is the Presentation
Axiom (PAx), studied, for instance, by Aczel \cite{Aczel} in the context of constructive set theory.
\end{example}

\begin{example}[Shulman]
$(\Top,\mathcal{O})$ satisfies WISC, using AC in $\Set$.
\end{example}

Choice may be more than is necessary here; it would be interesting to see if weaker 
choice principles in the site $(\Set,surjections)$ are enough to prove WISC for 
$(\Top,\mathcal{O})$ or other concrete sites.

\begin{lemma}
If $(S,J)$ satisfies WISC, then so does $(S,J_{un})$.
\end{lemma}

It is instructive to consider an example where WISC fails in a non-artificial 
way.  The category of sets and surjections with all arrows covers clearly 
doesn't satisfy WISC,  but is contrived and not a `useful' sort of category.
For the moment, assume the existence
of a Grothendieck universe $\mathbb{U}$ with cardinality $\lambda$,
 and let $\mathrm{Set}_\mathbb{U}$ refer to the category of
$\mathbb{U}$-small sets. Clearly we can define WISC relative to $\mathbb{U}$,
call it WISC${}_\mathbb{U}$.
 Let $G$ be a  $\mathbb{U}$-large group and $\mathbf{B}G$
the $\mathbb{U}$-large groupoid with one object associated to $G$. The boolean 
topos $\mathrm{Set}_\mathbb{U}^{\mathbf{B}G}$ of $\mathbb{U}$-small $G$-sets is 
a unary site with the class $epi$ of epimorphisms for covers. One could consider 
this topos as being an exotic sort of forcing construction.

\begin{proposition}\label{U-small_G-sets}
If $G$ has at least $\lambda$-many conjugacy classes of subgroups, 
then $(\mathrm{Set}_\mathbb{U}^{\mathbf{B}G},epi)$ does not satisfy 
WISC${}_\mathbb{U}$.
\end{proposition}
Alternatively, one could work in foundations where it is legitimate to discuss
a proper class-sized group, and then consider the topos of sets with an action
by this group. If there is a proper class of conjugacy classes of subgroups, then 
this topos with its canonical singleton pretopology will fail to satisfy WISC. Simple
examples of such groups are $\mathbb{Z}^\mathbb{U}$ (given a universe $\mathbb{U}$) and
$\mathbb{Z}^K$ (for some proper class $K$).

Recently, \cite{vdBerg_12} (relative to a large cardinal axiom) and \cite{Roberts_13}  (with no large cardinals)
have shown that the category of sets may fail to satisfy WISC. The models constructed
in \cite{Karaglia_12} are also conjectured to not satisfy WISC.

Perhaps of independent interest is a form of WISC with a bound: the weakly 
initial set for each category $J/A$ has cardinality less than some cardinal $\kappa$
(call this WISC${}_\kappa$).
Then one could consider, for example, sites where each 
object has a weakly initial finite or countable set of covers. Note that the condition
`enough projectives' is the case $\kappa = 2$.

\section{Internal categories}\label{internal_cats}

Internal categories were introduced in \cite{Ehresmann_63}, starting
with differentiable and topological categories (i.e. internal to $\Diff$ and 
$\Top$ respectively). We collect here the necessary definitions, terminology and 
notation. For a thorough recent account, see \cite{HDA5} or the encyclopedic 
\cite{Elephant}.

Fix a category $S$, referred to as the \emph{ambient category}.

\begin{definition}\label{def:cat}
An \emph{internal category} $X$ in a category $S$ is a diagram
\[
      X_1 \times_{X_0} X_1 \times_{X_0} X_1\rightrightarrows X_1 \times_{X_0} 
      X_1 \xrightarrow{m} X_1 \stackrel{s,t}{\st} X_0 \xrightarrow{e} X_1
\]
in $S$ such that the \emph{multiplication} $m$ is associative (we demand the 
limits in the 
diagram exist), the \emph{unit map} $e$
is a two-sided unit for $m$ and $s$ and $t$ are the usual \emph{source} and
\emph{target}.  An \emph{internal groupoid} is an internal category with an 
involution
\[
       (-)^{-1}\colon X_1 \to X_1
\]
satisfying the usual diagrams for an inverse.
\end{definition}

Since multiplication is associative, there is a well-defined map 
$X_1 \times_{X_0} X_1 \times_{X_0} X_1 \to X_1$, which will also be denoted by 
$m$. The pullback in the diagram in definition \ref{def:cat} is
\[
       \xymatrix{
               X_1 \times_{X_0} X_1 \ar[r] \ar[d] & X_1 \ar[d]^-{s}\\
               X_1 \ar[r]_-{t} & X_0\;.
       }
\]
and the double pullback is the limit of 
$X_1 \stackrel{t}{\rightarrow} X_0 \stackrel{s}{\leftarrow} X_1 
\stackrel{t}{\rightarrow} X_0 \stackrel{s}{\leftarrow}X_0$.
These, and pullbacks like these (where source is pulled back along target), will 
occur often. If confusion can arise, the maps in question will be explicity 
written, as in $X_1 \times_{s,X_0,t} X_1$. One usually sees the requirement that 
$S$ is finitely complete in order to define internal categories. This is not 
strictly necessary, and not true in the well-studied case of $S = \Diff$, the 
category of smooth manifolds.

Often an internal category will be denoted $X_1 \st X_0$, the arrows $m,s,t,e$ 
(and $(-)^{-1}$) will be referred to as \emph{structure maps} and $X_1$ and 
$X_0$ called the object of arrows and the object of objects respectively. For 
example, if $S = \Top$, we have the space of arrows and the space of objects, 
for $S = \Grp$ we have the group of arrows and so on.

\begin{example}\label{eg:cech_gpd}
If $X \to Y$ is an arrow in $S$ admitting iterated kernel pairs, there is an 
internal groupoid $\check{C}(X)$ with 
$\check{C}(X)_0 = X$, $\check{C}(X)_1 = X \times_Y X$,
source and target are projection on first and second factor, and the 
multiplication is projecting out the middle factor in $X \times_Y X \times_Y X$. 
This groupoid is called the \emph{\v Cech groupoid} of the map $X \to Y$. The origin of 
the name is that in $\Top$, for maps of the form $\coprod_I U_i \to Y$ (arising 
from an open cover), the \v Cech groupoid $\check{C}(\coprod_I U_i)$ appears in 
the definition of \v Cech cohomology.
\end{example}

\begin{example}\label{eg:disc-codisc_gpd}
Let $S$ be a category with binary products. For each object $A \in S$ there is an internal groupoid
$\disc(A)$ which has $\disc(A)_1 = \disc(A)_0 = A$ and all structure maps equal 
to $id_A$. Such a category is called \emph{discrete}. 
There is also an internal groupoid $\codisc(A)$ with 
\[
\codisc(A)_0 = A,\
\codisc(A)_1 = A \times A
\]
and where source and target are projections on the 
first and second factor respectively. Such a groupoid is called 
\emph{codiscrete}. 
\end{example}

\begin{definition}\label{def:functor}
Given internal categories $X$ and $Y$ in $S$, an \emph{internal functor}
\mbox{$f: X \to Y$} is a pair of maps
\[
f_0\colon X_0 \to Y_0 \quad\textrm{and}\quad f_1\colon X_1 \to Y_1
\]
called the object and arrow component respectively.  Both components are 
required to commute with all the structure maps.
\end{definition}

\begin{example}
If $A\to C$ and $B\to C$ are maps admitting iterated kernel pairs, and $A \to B$ 
is a map over $C$, there is a functor $\check{C}(A) \to \check{C}(B)$.
\end{example}

\begin{example}\label{functors2discrete}
If $(S,J)$ is a subcanonical unary site, and $U \to A$ is a cover, 
a functor $\check{C}(U) \to \disc(B)$
gives a unique arrow $A\to B$. This follows immediately from the fact $A$ is the 
colimit of the diagram underlying $\check{C}(U)$.
\end{example}

\begin{definition}\label{def:nat_iso}
Given internal categories $X,Y$ and internal functors $f,g\colon X \to Y$, an
\emph{internal natural transformation} (or simply \emph{transformation})
\[
       a\colon f \Rightarrow g
\]
is a map $a\colon X_0 \to Y_1$ such that $s \circ a = f_0,\ t\circ a = g_0$ and 
the following diagram commutes
\begin{equation}\label{diag:naturality}
       \xymatrix{
               X_1 \ar[r]^-{(g_1,a\circ s)} \ar[d]_{(a \circ t,f_1)} &
               Y_1 \times_{Y_0} Y_1 \ar[d]^{m} \\
               Y_1 \times_{Y_0} Y_1 \ar[r]^-{m} & Y_1
       }
\end{equation}
expressing the naturality of $a$. 
\end{definition}

Internal categories (resp. groupoids), functors and transformations 
in a locally small category $S$ form a 
locally small  2-category $\Cat(S)$ (resp. $\Gpd(S)$) \cite{Ehresmann_63}. There 
is clearly an inclusion 2-functor $\Gpd(S) \to \Cat(S)$.
Also, $\disc$ and $\codisc$, described in example \ref{eg:disc-codisc_gpd}, are 
2-functors $S \to \Gpd(S)$, whose underlying functors
are left and right adjoint to the functor
\[
\Obj\colon\Cat(S)_{\leq 1} \to S,\qquad (X_1\st X_0)\mapsto X_0.
\]
Here $\Cat(S)_{\leq 1}$ is the 1-category underlying the 2-category 
$\Cat(S)$. Hence for an internal category $X$ in $S$, there are functors 
$\disc(X_0) \to X$ and $X \to \codisc(X_0)$, the arrow component of the latter 
being $(s,t):X_1\to X_0^2$.

We say a natural transformation is a \emph{natural isomorphism} if it has an 
inverse with respect to vertical composition. Clearly there is no distinction 
between natural transformations and natural isomorphisms when the codomain of 
the functors is an internal groupoid. 
We can reformulate the naturality diagram (\ref{diag:naturality}) in the case 
that $a$ is a natural isomorphism. Denote by $-a$ the inverse of $a$.
Then the diagram (\ref{diag:naturality}) commutes if and only if the diagram 
\begin{equation}\label{naturality}
       \xymatrix{
             X_0 \times_{X_0} X_1 \times_{X_0} X_0 
             \ar[rr]^{-a\times f_1 \times a}
			 \ar[d]_{\simeq} &&Y_1 \times_{Y_0} Y_1 \times_{Y_0} Y_1 \ar[d]^m \\
             X_1 \ar[rr]_{g_1} && Y_1
       }
\end{equation}
commutes, a fact we will use several times.

\begin{example}
If $X$ is a category in $S$, $A$ is an object of $S$ and $f,g:X \to \codisc(A)$ 
are functors, there is a unique natural isomorphism 
$f\stackrel{\sim}{\Rightarrow} g$.
\end{example}

\begin{definition}
An \emph{internal} or \emph{strong equivalence} of internal categories is an
equivalence in the 2-category of internal categories. That is, an internal functor $f
\colon X\to Y$ such that there is a functor $f'\colon Y\to X$ and natural isomorphisms
$f\circ f' \Rightarrow \id_Y$, $f'\circ f \Rightarrow \id_X$.
\end{definition}

\begin{definition}\label{induced_cat}
For an internal category $X$ and a map $p:M\to X_0$ in $S$ the \emph{base change of $X$ along
$p$} is any category $X[M]$ with object of objects $M$ and object of arrows given by the pullback
\[
\xymatrix{
               M^2 \times_{X_0^2} X_1 \ar[r] \ar[d]  & X_1 \ar[d]^{(s,t)} \\
               M^2 \ar[r]_{p^2} & X_0^2
       }
\]
\end{definition}

 If $C\subset \Cat(S)$ denotes a full sub-2-category and if the base change
 along any map in a given class $K$ of maps exists in $C$ for all objects
of $C$, then we say $C$ \emph{admits base change along maps in $K$},
or simply \emph{admits base change for $K$}.

\begin{remark}
In all that follows, `category' will mean object of $C$ and similarly for 
`functor' and `natural transformation/isomorphism'.
\end{remark}

The strict pullback of internal categories
\[
       \xymatrix{
               X \times_Y Z \ar[r] \ar[d] & Z \ar[d] \\
               X \ar[r] & Y
       }
\]
when it exists, is the internal category with objects $X_0 \times_{Y_0} Z_0$, arrows
$X_1 \times_{Y_1} Z_1$, and all structure maps given componentwise by those of
$X$ and $Z$. Often we will be able to prove that certain pullbacks exist because of
conditions on various component maps in $S$. We do not assume that all strict
pullbacks of internal categories exists in our chosen $C$.

It follows immediately from definition \ref{induced_cat} that given maps $N\to M$ and $M\to X_0$,
there is a  canonical isomorphism
\begin{equation}\label{induced_cat_1}
       X[M][N] \simeq X[N].
\end{equation}
with object component the identity map, when these base changes exist.

\begin{remark}\label{remark:pullback_of_x_along_x}
If we agree to follow the convention that $M \times_N N = M$ is the pullback along
the identity arrow $\id_N$, then $X[X_0] = X$. This also simplifies other results of this
paper, so will be adopted from now on.
\end{remark}

One consequence of this assumption is that the iterated fibre product
\[
       M\times_M M \times_M \ldots \times_M M,
\]
bracketed in any order, is \emph{equal} to $M$. We cannot, however, equate two
bracketings of a general iterated fibred product; they are only canonically isomorphic.

\begin{lemma}\label{strict_pullbacks_triv_J_fibrations}
       Let $Y\to X$ be a functor in $S$ and $j_0\colon U \to X_0$ a map. 
       If the base change along $j_0$ exists, the following square is a strict pullback
       \[
       \xymatrix{
               Y[Y_0\times_{X_0}U] \ar[r] \ar[d] & X[U] \ar[d]^j \\
               Y \ar[r] & X
       }
       \]
       assuming it exists.
\end{lemma}
\proof
Since base change along $j_0$ exists, we know that we have the functor 
$Y[Y_0\times_{X_0}U] \to Y$, we just need to show it is a strict pullback of $j$. On
the level of objects this is clear, and on the level of arrows, we have
\begin{align*}
(Y_0\times_{X_0}U)^2 \times_{Y_0^2}Y_1 &\simeq U^2\times_{X_0^2} Y_1\\
	 &\simeq (U^2\times_{X_0^2}X_1) \times_{X_1}Y_1 \\
	 &\simeq X[U]_1\times_{X_1}Y_1
\end{align*}
so the square is a pullback.
\endofproof

We are interested in 2-categories $C$ which admits base change for a
given pretopology $J$ on $S$, which we shall cover in more detail in section \ref{examples}.

Equivalences in $\Cat$---assuming the axiom of choice---are precisely the fully
faithful, essentially surjective functors. For internal categories, however, this is not
the case. In addition, we need to make use of a pretopology to make the `surjective'
part of `essentially surjective' meaningful.

\begin{definition}\label{def:weak_equiv}
Let $(S,J)$ be a unary site. An internal functor $f:X \to Y$ in $S$ is called
\begin{enumerate}
       \item
       \emph{fully faithful} if
       \[
               \xymatrix{
                       X_1 \ar[r]^{f_1} \ar[d]_{(s,t)} & Y_1 \ar[d]^{(s,t)}\\
                       X_0 \times X_0 \ar[r]_{f_0 \times f_0} & Y_0 \times Y_0
               }
       \]
       is a pullback diagram;
       \item
       \emph{$J$-locally split} if there is a $J$-cover $U\to Y_0$ and a diagram
       \[
	\xymatrix{
		Y[U] \ar[d]_{\bar f} \ar@/^.5pc/[dr]_{\ }="s1"^{u}& \\
		X\ar[r]_{f}^(.33){\ }="t1"&Y 
		\ar@{=>}"s1";"t1"
	}
	\]
	commuting up to a natural isomorphism;
	\item
       a \emph{$J$-equivalence} if it is fully faithful and $J$-locally split.
\end{enumerate}
The class of $J$-equivalences will be denoted $W_J$.  If mention of $J$ is
suppressed, they will be called \emph{weak equivalences}. 
\end{definition}

\begin{remark}
There is another defintion of full faithfulness for internal categories, namely
that of a functor $f\colon Z\to Y$ being \emph{representably fully faithful}. 
This means that for all categories $Z$, the functor 
\[
	f_\ast\colon \Cat(S)(Z,X) \to \Cat(S)(Z,Y)
\]
is fully faithful. It is a well-known result that these two notions coincide, so we
shall use either characterisation as needed.
\end{remark}

\begin{lemma}
If $f:X \to Y$ is a fully faithful functor such that $f_0$ is in $J$, then $f$ is $J$-locally split. 
\end{lemma}

That is, the canonical functor $X[U] \to X$ is a $J$-equivalence whenever the
base change exists. Also, we do not require that $J$ is subcanonical. We record
here a useful lemma.

\begin{lemma}\label{rep_ff_functors_closed_under_iso}
Given a  fully faithful functor $f\colon X \to Y$ in $C$ and a natural
isomorphism $f \Rightarrow g$, the functor $g$ is also fully faithful. In 
particular, an internal equivalence is  fully faithful.
\end{lemma}
\proof
This is a simple application of the definition of representable full faithfulness
and the fact that the result is true in $\Cat$.
\endofproof

The first definition of weak equivalence of internal categories along the lines we are 
considering appeared in \cite{Bunge-Pare_79} for $S$
a regular category, and $J$ the class of regular epimorphisms (i.e.~$\can$), in the
context of stacks and indexed categories. This was later generalised in
\cite{Everaert_et_al_05} to more general finitely complete sites to discuss
model structures on the category of internal categories. Both work only 
with saturated singleton pretopologies.

Note that when $S$ is finitely complete, the object
$X_1^{iso} \into X_1$ of isomorphisms of a category $X$ can be constructed as a finite limit
\cite{Bunge-Pare_79}, and in the case when $X$ is a groupoid we have $X_1^{iso} \simeq X_1$.

\begin{definition}{\cite{Bunge-Pare_79,Everaert_et_al_05}}\label{def:PB_weak_equiv}
	For a finitely complete unary site $(S,J)$ with $J$ saturated, a functor $f$ is called
       \emph{essentially $J$-surjective} if the arrow labelled $\circledast$ below is in $J$.
       \[
               \xymatrix{
                       &\ar[dl] X_0 \times_{Y_0} Y_1^{iso} \ar@/^1pc/[ddr]^\circledast \ar[d]&\\
                       X_0 \ar[d]_{f_0} & \ar[dl]^s Y_1^{iso} \ar[dr]_t  &\\
                       Y_0 && Y_0
               }
       \]
       A functor is called a \emph{Bunge-Par\'e $J$-equivalence} if it is 
       fully faithful and essentially $J$-surjective. Denote the class of such
       maps by $W_J^{BP}$.
\end{definition}

Definition \ref{def:weak_equiv} is equivalent to the one in \cite{Bunge-Pare_79,Everaert_et_al_05} 
in the sites they consider but seems more appropriate for sites without all finite
limits. Also, definition \ref{def:weak_equiv} makes sense in 2-categories other than $\Cat(S)$
or sub-2-categories thereof.


\begin{proposition}\label{BP_equiv_iff_weak_equiv}
Let $(S,J)$ be a finitely complete unary site with $J$ saturated.
Then a functor is a $J$-equivalence if and only if it is a Bunge-Par\'e
$J$-equivalence.
\end{proposition}

\proof

Let $f\colon X \to Y$ be a Bunge-Par\'e $J$-equivalence, and consider the
$J$-cover given by the map $U := X_0 \times_{Y_0} Y_1^{iso} \to Y_0$. Denote by
$\iota\colon U\to Y_1^{iso}$ the projection on the second factor, by $-\iota$ the composite
of $\iota$ with the inversion map $(-)^{-1}$ and by $s_0\colon U\to X_0$
the projection on the first factor.
The arrow $s_0$ will be the object component of a functor $s\colon Y[U] \to X$, 
we need to define the arrow component $s_1$. Consider the composite
\begin{align*}
	Y[U]_1 \simeq U\times_{Y_0} Y_1 \times_{Y_0} U \xrightarrow{(s,\iota)\times\id\times(-\iota,s)}
		(X_0 \times_{Y_0} Y_1^{iso}) \times_{Y_0} Y_1 \times_{Y_0} ( Y_1^{iso} \times_{Y_0} X_0) \\
	\hookrightarrow X_0 \times_{Y_0} Y_3 \times_{Y_0} X_0 \xrightarrow{\id\times m\times\id} 
		X_0 \times_{Y_0} Y_1 \times_{Y_0} X_0 \simeq X_1
\end{align*}
where the last isomorphism arises from $f$ being fully faithful. It is clear that this
commutes with source and target, because these are given by projection on the first and last
factor at each step. To see that it respects identities and composition, one can use generalised
elements and the fact that the $\iota$ component will cancel with the 
$-\iota = (-)^{-1}\circ \iota$ component.

We define the natural isomorphism $f\circ s \Rightarrow j$ (here $j\colon Y[U] \to Y$ 
is the canonical functor) to have component $\iota$ as denoted above. 
Notice that the composite $f_1\circ s_1$ is just
\[
       Y[U]_1 \simeq U \times_{Y_0} Y_1 \times_{Y_0} U \xrightarrow{\iota\times\id
\times -\iota}  Y_1^{iso} \times_{Y_0} Y_1 \times_{Y_0}  Y_1^{iso}  \hookrightarrow
       Y_3  \xrightarrow{m} Y_1.
\]
Since the arrow component of $Y[U] \to Y$ is $U \times_{Y_0} Y_1 \times_{Y_0} U
\xrightarrow{\pr_2} Y_1$, $\iota$ is indeed a natural isomorphism using the diagram
(\ref{naturality}). Thus a Bunge-Par\'e $J$-equivalence is a $J$-equivalence.

In the other direction, given a $J$-equivalence $f\colon X\to Y$, we have a 
$J$-cover $j\colon U\to Y_0$ and a map $(\overline{f},a)\colon U \to X_0 \times Y_1^{iso}$ 
such that $j = (t\circ pr_2)\circ(\overline{f},a)$. Since $J$ is saturated, $(t\circ pr_2)\in J$ 
and hence $f$ is a Buge-Par\'e $J$-equivalence.
\endofproof

We can thus use definition \ref{def:weak_equiv} as we like, and it will still
refer to the same sorts of weak equivalences that appear in the literature.

\section{Anafunctors}\label{anafunctors}

We now let $J$ be a \emph{subcanonical} singleton pretopology on the ambient
category $S$. In this section we assume that $C\into \Cat(S)$ admits base change along
arrows in the given pretopology $J$. This is a slight generalisation of what is
considered in \cite{Bartels}, where only $C = \Cat(S)$ is considered.

\begin{definition}{\cite{Makkai,Bartels}}\label{def:anafunctor}
An \emph{anafunctor} in $(S,J)$ from a category $X$ to a category $Y$ consists of a
$J$-cover $(U \to X_0)$ and an internal functor
\[
       f\colon X[U] \to Y.
\]
Since $X[U]$ is an object of $C$, an anafunctor is a span in $C$, and
can be denoted
\[
       (U,f)\colon X \gento Y.
\]
\end{definition}

\begin{example}\label{eg:ordinary_functor}
For an internal functor $f\colon X \to Y$ in $S$, define the anafunctor $(X_0,f)\colon
X \gento Y$ as the following span
\[
       X \xleftarrow{=} X[X_0]  \xrightarrow{f} Y.
\]
We will blur the distinction between these two descriptions. If $f=id\colon X \to X$,
then $(X_0,id)$ will be denoted simply by $id_X$.
\end{example}

\begin{example}
If $U \to A$ is a cover in $(S,J)$ and $\mathbf{B}G$ is a groupoid with one object in
$S$ (i.e. a group in $S$), an anafunctor $(U,g)\colon\disc(A) \gento \mathbf{B}G$ is the
same thing as a \v Cech cocycle.
\end{example}

\begin{definition}{\cite{Makkai,Bartels}}
Let $(S,J)$ be a site and let
\[
(U,f),(V,g)\colon X \gento Y
\]
be anafunctors in $S$. A \emph{transformation}
\[
       \alpha\colon (U,f) \Rightarrow (V,g)
\]
from $(U,f)$ to $(V,g)$ is a natural transformation
\[
       \xymatrix{
               & \ar[dl] X[U\times_{X_0}V]  \ar[dr] & \\
               X[U] \ar[dr]_f & \stackrel{\alpha}{\Rightarrow} & X[V] \ar[dl]^g\\
               & Y &
       }
\]
If $\alpha$ is a natural isomorphism, then $\alpha$ will be called an
\emph{isotransformation}. In that case we say $(U,f)$ is
isomorphic to $(V,g)$. Clearly all transformations between anafunctors between
internal groupoids are isotransformations.
\end{definition}

\begin{example}\label{eg:ordinary_transf}
Given functors $f,g\colon X \to Y$ between categories in $S$, and a natural
transformation $a\colon f \Rightarrow g$, there is a transformation $a\colon (X_0,f)
\Rightarrow (X_0,g)$ of anafunctors, given by the component $X_0\times_{X_0}X_0
= X_0 \xrightarrow{a} Y_1$.
\end{example}

\begin{example}
If $(U,g),(V,h)\colon \disc(A) \gento \mathbf{B}G$ are two \v Cech cocycles, a
transformation between them is a coboundary on the cover $U\times_A V\to A$.
\end{example}

\begin{example}
Let $(U,f)\colon X \gento Y$ be an anafunctor in $S$. There is an isotransformation
$1_{(U,f)}\colon (U,f) \Rightarrow (U,f)$ called the \emph{identity transformation},
given by the natural transformation with component
\begin{equation}\label{id_transf_component}
       U \times_{X_0} U \simeq (U \times U) \times_{X_0^2} X_0 \xrightarrow{id_U^2
\times e} X[U]_1 \xrightarrow{f_1} Y_1
\end{equation}
\end{example}

\begin{example}{\cite{Makkai}}\label{renaming_transf}
Given anafunctors $(U,f)\colon X\to Y$ and $(V,f\circ k)\colon X \to Y$ where $k
\colon V\to U$ is a cover (over $X_0$), a \emph{renaming transformation}
\[
(U,f)\Rightarrow(V,f\circ k)
\]
 is an isotransformation with component
\[
       1_{(U,f)}\circ (k\times \id):V\times_{X_0} U \to U\times_{X_0} U \to Y_1.
\]
(We also call its inverse for vertical composition a renaming transformation.) 
If $k$ is an isomorphism, then it will itself be referred to as a
\emph{renaming isomorphism}.
\end{example}

We define (following \cite{Bartels}) the composition of anafunctors as follows. Let
\[
(U,f)\colon X \gento Y \quad \textrm{and} \quad (V,g)\colon Y \gento Z
\]
be anafunctors in the site $(S,J)$. Their composite $(V,g)\circ(U,f)$ is the composite
span defined in the usual way. It is again a span in $C$:
\[
       \xymatrix{
               && \ar[dl] X[U\times_{Y_0}V]  \ar[dr]^{f^V} & \\
               &\ar[dl]X[U] \ar[dr]_f &  & Y[V] \ar[dl] \ar[dr]^g\\
               X&& Y &&Z
       }
\]
The square is a pullback by lemma \ref{strict_pullbacks_triv_J_fibrations}
(which exists because $V\to
Y_0$ is a cover), and the resulting span is an anafunctor because $V \to Y_0$, 
hence $U\times_{Y_0}V\to X_0$, are covers, and using the isomorphism (\ref{induced_cat_1}). We will
sometimes denote the composite by $(U\times_{Y_0}V,g\circ f^V)$.

Here we are using the fact we have specified pullbacks of covers in $S$. Without this
we would not end up with a bicategory (see theorem \ref{anafunctors_are_a_bicat}), 
but what \cite{Makkai} calls an
\emph{anabicategory}. This is similar to a bicategory, but composition
and other structural maps are only anafunctors, not functors.

Consider the special case when $V = Y_0$, so that $(Y_0,g)$ is just an ordinary
functor. Then there is a renaming transformation (the identity transformation!)
$(Y_0,g)\circ(U,f) \Rightarrow (U,g\circ f)$, using the equality $U \times_{Y_0} Y_0= U$
(by remark \ref{remark:pullback_of_x_along_x}).
 If we let $g=\id_Y$, then we see that $(Y_0,\id_Y)$ is a strict unit on the left for
anafunctor composition. Similarly, considering $(V,g)\circ(Y_0,\id)$,
we see that $(Y_0,\id_Y)$ is a two-sided strict unit for anafunctor composition. In
fact, we have also proved
\begin{lemma}\label{coherent_composition}
Given two functors $f\colon X\to Y$, $g\colon Y \to Z$ in $S$, their composition as
anafunctors is equal to their composition as functors:
\[
       (Y_0,g)\circ(X_0,f) = (X_0,g\circ f).
\]
\end{lemma}

As a concrete and relevant example of a renaming transformation we can consider
the triple composition of anafunctors
\begin{align*}
       (U,f)\colon & X \gento Y,\\
       (V,g)\colon & Y \gento Z,\\
       (W,h)\colon & Z \gento A.
\end{align*}
The two possibilities of composing these are
\[
       \left((U\times_{Y_0} V)\times_{Z_0}W,h\circ(gf^V)^W\right)\quad \text{and}\quad \left(U
\times_{Y_0} (V\times_{Z_0} W),h\circ g^W\circ f^{V\times_{Z_0}W}\right).
\]
\begin{lemma}\label{lemma:associator}
The unique isomorphism $(U\times_{Y_0} V)\times_{Z_0}W \simeq U\times_{Y_0} (V
\times_{Z_0} W)$ commuting with the various projections is a
renaming isomorphism. The isotransformation arising from this renaming
transformation is called the \emph{associator}.
\end{lemma}

A simple but useful criterion for describing isotransformations where one of the
anafunctors involved is a functor is as follows.

\begin{lemma}\label{anafun_iso2_fun}
An anafunctor $(V,g)\colon X \gento Y$ is isomorphic to a functor $(X_0,f)\colon X \gento Y$ if
and only if there is a natural isomorphism
\[
       \xymatrix{
               & \ar[dl] X[V] \ar[dr]^g \\
               X \ar@/_1.5pc/[rr]_(.6){f}& \stackrel{\sim}{\Rightarrow} & Y
       }
\]
\end{lemma}

Just as there is a vertical composition of natural transformations between internal functors,
there is a vertical composition of transformations between internal anafunctors \cite{Bartels}.
This is where the subcanonicity of $J$ will be used in order to construct a map
locally over some cover. Consider the following diagram
\[
       \xymatrix{
               && \ar[dl] X[U\times_{X_0} V\times_{X_0} W]  \ar[dr]\\
               &        \ar[dl] X[U\times_{X_0} V] \ar[dr]   &
               & \ar[dl] X[V\times_{X_0} W] \ar[dr]\\
               X[U] \ar[drr]_f & \stackrel{a}{\Rightarrow} & X[V] \ar[d]^g&  \stackrel{b}
{\Rightarrow} & X[W] \ar[dll]^h \\
               &&Y&&
       }
\]
We can form a natural transformation between the leftmost and the
rightmost composites as functors in $S$. This will have as its component the arrow
\[
\widetilde{ba}\colon U\times_{X_0} V\times_{X_0} W \xrightarrow{\id\times \Delta
\times \id}  U\times_{X_0}V\times_{X_0}V\times_{X_0} W \xrightarrow{a\times b}
Y_1\times_{Y_0}  Y_1 \xrightarrow{m} Y_1
\]
in $S$. Notice that the \v Cech groupoid of the cover
\begin{equation}\label{iterated_cover}
U\times_{X_0} V\times_{X_0} W \to U \times_{X_0} W
\end{equation}
is
\[
U\times_{X_0} V\times_{X_0} V\times_{X_0} W \st U\times_{X_0} V\times_{X_0} W,
\]
with source and target arising from the two projections $V\times_{X_0} V \to V$.
Denote this pair of parallel arrows
by $s,t\colon UV^2W \st UVW$ for brevity. In \cite{Bartels}, section 2.2.3, we find the
commuting diagram
\begin{equation}\label{tobys_diag}
       \xymatrix{
               UV^2W \ar[r]^t \ar[d]_s & UVW \ar[d]^{\widetilde{ba}}\\
               UVW \ar[r]_{\widetilde{ba}} & Y_1
       }
\end{equation}
(this can be checked by using generalised elements) and so we have a functor 
\[
\check{C}(U\times_{X_0} V\times_{X_0} W) \to \disc(Y_1).
\]
 Our pretopology $J$ is assumed to
be subcanonical, so example \ref{functors2discrete} gives us a unique arrow
$ba\colon U\times_{X_0} W \to Y_1$, which is the data for the composite of $a$ and $b$.

\begin{remark}
In the special case that $U\times_{X_0} V\times_{X_0} W \to U \times_{X_0} W$ is 
split (e.g. is an isomorphism), the composite transformation has
\[
U \times_{X_0} W\to U\times_{X_0} V\times_{X_0} W \xrightarrow{\widetilde{ba}} Y_1
\]
as its component arrow. In particular, this is the case if one of $a$ or $b$ is a
renaming transformation.
\end{remark}

\begin{example}\label{eg:transf_compose2}
Let $(U,f):X\gento Y$ be an anafunctor and $U'' \xrightarrow{j'} U' \xrightarrow{j} U$
successive refinements of $U \to X_0$ (e.g isomorphisms). Let $(U',f_{U'})$ and
$(U'',f_{U''})$ denote the composites of $f$ with $X[U'] \to X[U]$ and $X[U''] \to X[U]$
respectively. The arrow
\[
       U \times_{X_0} U'' \xrightarrow{\id_U\times j\circ j'} U \times_{X_0} U \to Y_1
\]
is the component for the composition of the isotransformations $(U,f)
\Rightarrow(U',f_{U'}),\Rightarrow(U'',f_{U''})$ described in example
\ref{renaming_transf}. Thus we can see that the composite of renaming
transformations associated to isomorphisms $\phi_1,\phi_2$ is simply the renaming
transformation associated to their composite $\phi_1\circ \phi_2$.
\end{example}

This can be used to show that the associator satisfies the necessary coherence conditions.

\begin{example}\label{eg:transf_compose1}
If $a\colon f\Rightarrow g,\ b\colon g\Rightarrow h$ are natural transformations
between functors $f,g,h\colon X\to Y$ in $S$, their composite as transformations
between anafunctors
\[
(X_0,f),(X_0,g),(X_0,h)\colon X\gento Y.
\]
is just their composite as natural transformations. This uses the equality
\[
X_0\times_{X_0} X_0\times_{X_0} X_0= X_0\times_{X_0} X_0 = X_0,
\]
which is due to our choice in remark \ref{remark:pullback_of_x_along_x} of canonical
pullbacks.
\end{example}

Even though we don't have pseudoinverses for weak equivalences of internal
categories, one might guess that the local splitting guaranteed to exist by definition
is actually more than just a splitting of sorts. This is in fact the case, if we use anafunctors.

\begin{lemma}\label{anafunctors-r-inverses}
Let $f\colon X \to Y$ be a $J$-equivalence in $S$. There is an anafunctor
\[
       (U,\bar{f})\colon Y \gento X
\]
and isotransformations
\begin{align*}
\iota\colon (X_0,f)\circ (U,\bar{f}) & \Rightarrow id_Y\\
\epsilon\colon(U,\bar{f})\circ (X_0,f) & \Rightarrow id_X
\end{align*}
\end{lemma}
\proof
We have the anafunctor $(U,\bar{f})$ by definition as $f$ is $J$-locally split.
Since the anafunctors $\id_X,\ \id_Y$ are actually
functors, we can use lemma \ref{anafun_iso2_fun}. Using the special case of
anafunctor composition when the second is a functor, this tells us that $\iota$ will be
given by a natural isomorphism
\[
       \xymatrix{
               & X \ar[dr]^{f}_(0.2){\ }="s" & \\
               Y[U] \ar[rr]^{\ }="t" \ar[ur]^{\bar{f}} && Y
                \ar@{=>}"s";"t"
       }
\]
with component $\iota\colon U \to Y_1$. Notice that the composite $f_1\circ \bar{f}_1$ is just
\[
       Y[U]_1 \simeq U \times_{Y_0} Y_1 \times_{Y_0} U \xrightarrow{\iota\times\id
\times -\iota}  Y_1 \times_{Y_0} Y_1 \times_{Y_0}  Y_1  \hookrightarrow
       Y_3  \xrightarrow{m} Y_1.
\]
Since the arrow component of $Y[U] \to Y$ is $U \times_{Y_0} Y_1 \times_{Y_0} U
\xrightarrow{\pr_2} Y_1$, $\iota$ is indeed a natural isomorphism using the diagram
(\ref{naturality}).

The other isotransformation $\epsilon$ is between $(X_0\times_{Y_0} U,\bar{f}\circ \pr_2)$ and
$(X_0,\id_X)$, and is given by the component
\[
\epsilon\colon X_0 \times_{X_0} X_0\times_{Y_0} U = X_0\times_{Y_0} U
\xrightarrow{\id\times (\bar{f}_0,\iota)} X_0\times_{Y_0} (X_0\times_{Y_0} Y_1) \simeq X_0^2
\times_{Y_0^2} Y_1 \simeq X_1
\]
The diagram
\[
       \xymatrix{
               (X_0\times_{Y_0^2} U)^2 \times_{X_0^2} X_1 \ar[d]_\simeq \ar[rr]^{\pr_2}
& &X_1 \ar[dd]^\simeq\\
               U \times_{Y_0} X_1 \times_{Y_0}U \ar[d]_{-\iota\times f\times\iota} & \\
               (X_0 \times_{Y_0} Y_1) \times_{Y_0} Y_1 \times_{Y_0} (Y_1 \times_{Y_0} X_0) 
               \ar[rr]_(.6){\id\times m \times \id}
               && X_0\times_{Y_0} Y_1 \times_{Y_0} X_0
       }
\]
commutes (a fact which can be checked using generalised elements), and using (\ref{naturality})
we see that $\epsilon$ is natural.
\endofproof

The first half of the following theorem is proposition 12 in \cite{Bartels}, and the
second half follows because all the constructions of categories involved in dealing
with anafunctors outlined above are still objects of $C$.

\begin{theorem}{\cite{Bartels}}\label{anafunctors_are_a_bicat}
For a site $(S,J)$ where $J$ is a subcanonical singleton pretopology, internal
categories, anafunctors and transformations form a bicategory $\Cat_\ana(S,J)$. If
we restrict attention to a full sub-2-category $C$ which admits base change for
arrows in $J$, we have an analogous full sub-bicategory $C_\ana(J)$.
\end{theorem}

In fact the bicategory $C_{ana}(J)$ fails to be a strict 2-category
only in the sense that the associator
is given by the non-identity isotransformation from lemma \ref{lemma:associator}.
All the other structure is strict.

There is a strict 2-functor $C_\ana(J) \to \Cat_\ana(S,J)$ which is an inclusion on
objects and fully faithful in the strictest sense, namely
being the identity functor on hom-categories. The following is the main
result of this section, and allows us to relate anafunctors to the localisations
considered in the next section.

\begin{proposition}\label{W-inverting_alpha}
There is a strict, identity-on-objects 2-functor
\[
\alpha_J\colon C \to C_\ana(J)
\]
sending $J$-equivalences to equivalences, and commuting with the respective
inclusions into $\Cat(S)$ and $\Cat_\ana(S,J)$.
\end{proposition}

\proof 
We define $\alpha_J$ to be the identity on objects, and as described in examples
\ref{eg:ordinary_functor}, \ref{eg:ordinary_transf} on 1-arrows and 2-arrows (i.e. functors
and transformations). We need first to show that this gives a functor $C(X,Y)
\to C_\ana(J)(X,Y)$. This is precisely the content of example
\ref{eg:transf_compose1}. Since the identity 1-cell on a category $X$ in
$C_\ana(J)$ is the image of the identity functor on $S$ in $C$, $\alpha_J$
respects identity 1-cells. Also, lemma \ref{coherent_composition} tells us that $
\alpha_J$ respects composition. That $\alpha_J$ sends $J$-equivalences to
equivalences is the content of lemma \ref{anafunctors-r-inverses}.
\endofproof

The 2-category $C$ is locally small (i.e.~enriched in small categories)
if $S$ itself is locally small (i.e.~enriched in sets),
but \emph{a priori} the collection of anafunctors $X\gento Y$ do not constitute 
a set for $S$ a large category.

\begin{proposition}
Let $(S,J)$ be a locally small, subcanonical unary site satisfying 
WISC and let $C$ admit base change along arrows in $J$. Then 
$C_\ana(J)$ is locally essentially small.
\end{proposition}

\proof
Given an object $A$ of $S$, let $I(A)$ be a weakly initial set for $J/A$. Consider the locally full 
sub-2-category of $C_\ana(J)$ with the same objects, and arrows those 
anafunctors $(U,f):X \gento Y$ such that $U \to X_0$ is in $I(X_0)$. Every 
anafunctor is then isomorphic, by example 
\ref{renaming_transf}, to one in this sub-2-category. The collection of anafunctors $(U,f):X \gento Y$
for a fixed $U$ forms a set, by local smallness of $C$, and similarly the collection
of transformations between a pair of anafunctors forms a set by local smallness of $S$.
\endofproof

Examples of locally small sites $(S,J)$ where $C_\ana(J)$ is not known to be locally 
essentially small are the category of sets from the model of ZF used in \cite{vdBerg_12},
the model of ZF constructed in \cite{Roberts_13} and the topos from 
proposition \ref{U-small_G-sets}. We note that local essential
smallness of $C_\ana(J)$ seems to be a condition just slightly weaker
than WISC.

\section{Localising bicategories at a class of 1-cells}\label{localisation}

Ultimately we are interesting in inverting all $J$-equivalences in $C$ and so
need to discuss what it means to add the formal pseudoinverses to a class of 1-cells
in a 2-category -- a process known as \emph{localisation}. This was done in
\cite{Pronk_96} for the more general case of a class of 1-cells in a bicategory, where
the resulting bicategory is constructed and its universal properties examined. 
The application in \emph{loc.~cit.}~is to show the
equivalence of various bicategories of stacks to localisations of 2-categories of
smooth, topological and algebraic groupoids.  The results of this article can be seen
as one-half of a generalisation of these results to more general sites.

\begin{definition}{\cite{Pronk_96}}
Let $B$ be a bicategory and $W \subset B_1$ a class of 1-cells. A \emph{localisation
of $B$ with respect to $W$} is a bicategory $B[W^{-1}]$ and a weak 2-functor
\[
       U \colon B \to B[W^{-1}]
\]
such that $U$ sends elements of $W$ to equivalences, and is universal with this
property i.e.~precomposition with $U$ gives an equivalence of bicategories
\[
       U^* \colon Hom(B[W^{-1}],D) \to Hom_W(B,D),
\]
where $Hom_W$ denotes the sub-bicategory of weak 2-functors that send elements
of $W$ to equivalences (call these \emph{$W$-inverting}, abusing notation slightly).
\end{definition}

The universal property means that $W$-inverting weak 2-functors $F\colon B \to D$
factor, up to an equivalence, through $B[W^{-1}]$, inducing an essentially unique
weak 2-functor  $\widetilde{F}\colon B[W^{-1}] \to D$.

\begin{definition}{\cite{Pronk_96}}\label{bicat_fracs}
Let $B$ be a  bicategory with a class $W$ of 1-cells. $W$ is said to \emph{admit
a right calculus of fractions} if it satisfies the following conditions
\begin{enumerate}
       \item[2CF1.]
               $W$ contains all equivalences
       \item[2CF2.]
               a) $W$ is closed under composition\\
               b) If $a\in W$ and there is an isomorphism $a \stackrel{\sim}{\Rightarrow} b$
               then    $b\in W$
       \item[2CF3.]
               For all $w\colon A' \to A,\ f\colon C \to A$ with $w\in W$ there exists a
               2-commutative square
               \[
               \xymatrix{
                       P \ar[dd]^v \ar[rr]^g && A'\ar[dd]^w_{\ }="s" \\
                       \\
                       C \ar[rr]^{f}="t" & & A
                       \ar@{=>}_{\simeq} "s"; "t"
               }
               \]
               with $v\in W$.
       \item[2CF4.]
               If $\alpha\colon w \circ f \Rightarrow w \circ g$ is a 2-arrow and $w\in W$
there is a 1-cell $v \in W$ and a 2-arrow $\beta\colon f\circ v \Rightarrow g \circ v$ such
that $\alpha\circ v = w \circ \beta$. Moreover: when $\alpha$ is an isomorphism, we
require $\beta$ to be an isomorphism too; when $v'$ and $\beta'$ form another such
pair, there exist 1-cells $u,\,u'$ such that $v\circ u$ and $v'\circ u'$ are in $W$, and
an isomorphism $\epsilon\colon v\circ u \Rightarrow v' \circ u'$ such that the following
diagram commutes:
               \begin{equation}\label{2cf4.diag}
                       \xymatrix{
                               f \circ v \circ u \ar@{=>}[rr]^{\beta\circ u}
                               \ar@{=>}[dd]_{f\circ \epsilon}^\simeq &&
                               g\circ v \circ u \ar@{=>}[dd]^{g\circ \epsilon}_\simeq \\
                               \\
                               f\circ v' \circ u' \ar@{=>}[rr]_{\beta'\circ u'} && g\circ v' \circ u'
                       }
               \end{equation}
\end{enumerate}
\end{definition}

For a bicategory $B$ with a calculus of right fractions, \cite{Pronk_96} 
constructs a localisation of $B$ as a bicategory of fractions; the 1-arrows are 
spans and the 2-arrows are  equivalence classes of bicategorical spans-of-spans diagrams.

From now on we shall refer to a calculus of right fractions as simply a calculus 
of fractions, and the resulting localisation constructed by Pronk as a 
bicategory of fractions. Since $B[W^{-1}]$ is defined only up to equivalence, it 
is of great interest to know when a bicategory $D$, in which elements of $W$ are 
sent to equivalences by a 2-functor $B \to D$, is equivalent to 
$B[W^{-1}]$. In particular, one might be interested in finding such an 
equivalent bicategory with a simpler description than that which appears in 
\cite{Pronk_96}.

\begin{proposition}{\cite{Pronk_96}}\label{comparison_thm}
A weak 2-functor $F:B \to D$ which sends elements of $W$ to equivalences induces
an equivalence of bicategories
\[
       \widetilde{F} \colon B[W^{-1}] \xrightarrow{\sim} D
\]
if the following conditions hold
\begin{enumerate}
       \item[EF1.]
               $F$ is essentially surjective,
       \item[EF2.]
               For every 1-cell $f \in D_1$ there are 1-cells $w \in W$ and $g\in B_1$
such that
               $Fg \stackrel{\sim}{\Rightarrow} f \circ Fw$,
       \item[EF3.]
               $F$ is locally fully faithful.
\end{enumerate}
\end{proposition}

Thanks are due to Matthieu Dupont for pointing out (in personal communication) 
that proposition \ref{comparison_thm} actually only holds in the one direction, not 
in both, as claimed in \emph{loc.~cit.}

The following is useful in showing a weak 2-functor sends weak equivalences to
equivalences, because this condition only needs to be checked on a class that is in
some sense cofinal in the weak equivalences.

\begin{proposition}\label{inverting_special_we}
Let $V \subset W$ be two classes of 1-cells in a bicategory $B$ such that for all
$w\in W$, there exists $v\in V$ and $s\in W$ and an invertible 2-cell
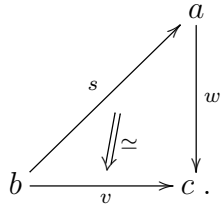
\[
       \xymatrix{
               && a \ar[dd]^w \\
               &  & \\
               b \ar[rr]_v^{\ }="t1" \ar[uurr]^s_{\ }="s1" && c\; .
               \ar@{=>}"s1";"t1"^{\simeq}
       }
\]
Then a weak 2-functor $F\colon B \to D$ that sends elements of $V$ to equivalences
also sends elements of $W$ to equivalences.
\end{proposition}
\proof
In the following the coherence arrows will be present, but unlabelled. It is 
enough to prove  that if in a bicategory $D$ with a class of maps $M$ (in our case 
$M=F(W)$) such that for all $w\in M$ there is an equivalence $v$ and an 
isomorphism $\alpha$,
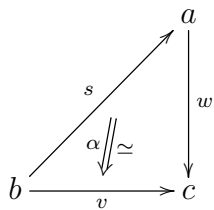
\[
       \xymatrix{
               && a \ar[dd]^w \\
               &  & \\
               b \ar[rr]_v^{\ }="t1" \ar[uurr]^s_{\ }="s1" && c
               \ar@{=>}"s1";"t1"^{\simeq}_\alpha
       }
\]
where $s\in M$, then all elements of $M$ are also equivalences. 

Let $\bar v$ be a pseudoinverse 
for $v$ and let $j = s \circ \bar v$. Then there is sequence of isomorphisms
\[
       w\circ j \Rightarrow (w\circ s)\circ \bar v \Rightarrow v \circ \bar v 
       \Rightarrow I.
\]

Since $s\in M$, there is an  equivalence $u$, $t\in M$ and an isomorphism 
$\beta$ giving the following diagram
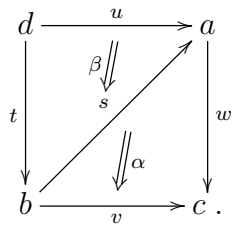
\[
       \xymatrix{
               d \ar[dd]_{t} \ar[rr]^{u}_{\ }="s2" && a \ar[dd]^w \\
               & & \\
               b \ar[rr]_v^{\ }="t1" \ar[uurr]^s="t2"_{\ }="s1" && c \; .
               \ar@{=>}"s1";"t1"^\alpha
               \ar@{=>}"s2";"t2"_\beta
       }
\]
Let $\bar u$ be a pseudoinverse of $u$. We know from the first part of the proof 
that we have a pseudosection $k = t\circ \bar u$ of $s$, with an isomorphism 
$s \circ k \Rightarrow I$. We then have the following sequence of isomorphisms:
\[
	j\circ w 
	= (s\circ \bar v) \circ w 
	\Rightarrow ((s\circ \bar v) \circ w) \circ (s \circ k)
	\Rightarrow s \circ ((\bar v \circ v) \circ (t\circ \bar u)) 
	\Rightarrow (s\circ t) \circ u
	\Rightarrow \bar u \circ u 
	\Rightarrow I.
\]
Thus all elements of $M$ are equivalences.
\endofproof

\section{2-categories of internal categories admit bicategories of fractions}\label{main}

In this section we prove the result that $C\into \Cat(S)$ admits a calculus of 
fractions for the $J$-equivalences, where $J$ is a 
singleton pretopology on $S$.

The following is the first main theorem of the paper, and subsumes a number of other, 
similar theorems throughout the literature (see section \ref{examples} for details).

\begin{theorem}\label{bicat_frac_exists}
Let $S$ be a category with a singleton pretopology $J$. Assume the full sub-2-category 
$C \into \Cat(S)$ admits base change along maps in $J$.  Then $C$ admits a 
right calculus of fractions for the class $W_J$ of $J$-equivalences.
\end{theorem}
\proof We show the conditions of definition \ref{bicat_fracs} hold.
\begin{itemize}
\item[2CF1.]  An internal equivalence is clearly $J$-locally split. Lemma 
\ref{rep_ff_functors_closed_under_iso} gives us the rest.

\item[2CF2.] \begin{itemize}
	\item[a)] That the composition of fully faithful functors is again fully faithful is trivial.
Consider the composition $g\circ f$ of two $J$-locally split functors,
\[
	\xymatrix{
		Y[U] \ar[d] \ar@/^.5pc/[dr]_{\ }="s1"^{u}&Z[V] \ar[d]\ar@/^.5pc/[dr]_(.5){\ }="s2"^{v}& \\
		X\ar[r]_{f}^(.33){\ }="t1"&Y \ar[r]_{g}^(.33){\ }="t2" & Z
		\ar@{=>}"s1";"t1"
		\ar@{=>}"s2";"t2"
	}
\]
By lemma \ref{strict_pullbacks_triv_J_fibrations} the functor $u$ pulls back to a functor
$Z[U\times_{Y_0}V] \to Z[V]$. The composite $Z[U\times_{Y_0}V] \to Z$ 
is  fully faithful with object component in $J$, hence $g\circ f$ is $J$-locally split.

\item[b)] Lemma  \ref{rep_ff_functors_closed_under_iso} tells us that fully faithful functors 
are closed under isomorphism, so we just need to show $J$-locally split functors are closed 
under isomorphism.

Let $w,f\colon X\to Y$ be functors and $a\colon w \Rightarrow f$ be a natural
isomorphism. First, let $w$ be $J$-locally split. It is immediate from the diagram
\[
	\xymatrix{
		Y[U] \ar[dd] \ar@/^.7pc/[ddrr]_{\ }="s1"^{u} \\
		\\
		X\ar@/^1pc/[rr]^{w}="t1"_{\ }="s2" \ar@/_1pc/[rr]_{f}^{\ }="t2"
		&&Y
		\ar@{=>}"s1";"t1"
		\ar@{=>}"s2";"t2"^{a}
	}
\]
that $f$ is also $J$-locally split.
\end{itemize}

\item[2CF3.] Let $w\colon X\to Y$ be a $J$-equivalence, and let $f\colon Z\to Y$ be a functor. From the definition of $J$-locally split, we have the diagram
\[
	\xymatrix{
		Y[U] \ar[d] \ar@/^.5pc/[dr]_{\ }="s1"^{u}& \\
		X\ar[r]_{w}^(.33){\ }="t1"&Y 
		\ar@{=>}"s1";"t1"
	}
\]
We can use lemma \ref{strict_pullbacks_triv_J_fibrations} to pull $u$ back along $f$ to get
a 2-commuting diagram
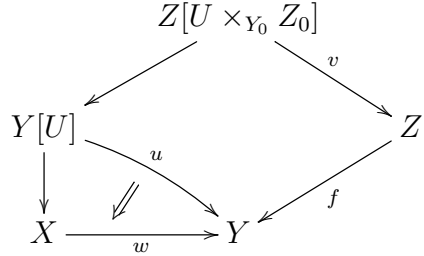
\[
	\xymatrix{
		& Z[U\times_{Y_0} Z_0] \ar[dr]^{v} \ar[dl] \\
		Y[U] \ar[d] \ar@/^.5pc/[dr]_{\ }="s1"^{u}& &Z \ar[dl]^f\\
		X\ar[r]_{w}^(.33){\ }="t1"&Y 
		\ar@{=>}"s1";"t1"
	}
\]
with $v\in W_J$ as required.

\item[2CF4.] 
Since $J$-equivalences are representably fully faithful, given
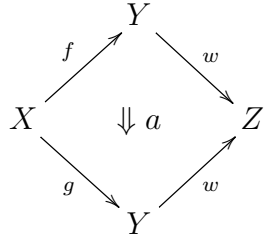
\[
       \xymatrix{
               &Y \ar[dr]^w \\
               X \ar[ur]^f \ar[dr]_g & \Downarrow a & Z\\
               & Y \ar[ur]_w
       }
\]
where $w\in W_J$, there is a unique $a'\colon f
\Rightarrow g$ such that
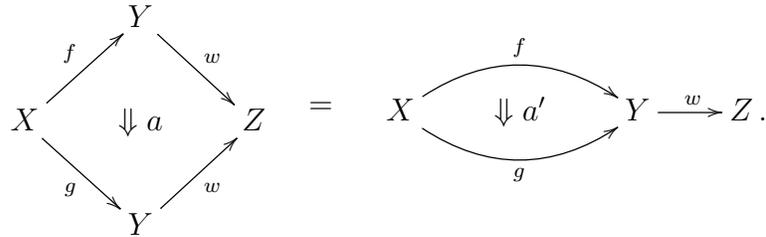
\[
        \raisebox{36pt}{
        \xymatrix{
               &Y \ar[dr]^w \\
               X \ar[ur]^f \ar[dr]_g & \Downarrow a & Z\\
               & Y \ar[ur]_w
       }
       }
       \equals
        \raisebox{36pt}{
       \xymatrix{
       &&\\
       X \ar@/^1.5pc/[rr]^f \ar@/_1.5pc/[rr]_g&\Downarrow a'& Y \ar[r]^w & Z
\,.
       }
       }
\]
The existence of $a'$ is the first half of 2CF4, where $v=\id_X$. Note that if $a$ is 
an isomorphism, so if $a'$, since $w$ is representably fully faithful.
Given $v'\colon W\to X \in W_J$ such that
there is a transformation
\[
       \xymatrix{
               &X \ar[dr]^f \\
               W \ar[ur]^{v'} \ar[dr]_{v'} & \Downarrow b & Y\\
               & X \ar[ur]_g
       }
\]
satisfying
\begin{align}\label{antiwhisker_eqn}
        \raisebox{36pt}{
       \xymatrix{
               &X \ar[dr]^f \\
               W \ar[ur]^{v'} \ar[dr]_{v'} & \Downarrow b & Y \ar[r]^w & Z\\
               & X \ar[ur]_g
       }
       }
       \equals &
        \raisebox{36pt}{
        \xymatrix{
               &&Y \ar[dr]^w \\
               W \ar[r]^{v'} &X \ar[ur]^f \ar[dr]_g & \Downarrow a & Z\\
               && Y \ar[ur]_w
       }
        }      \nonumber \\
        \equals &
        \raisebox{36pt}{
        \xymatrix{
               &&\\
               W\ar[r]^{v'}&X \ar@/^1.5pc/[rr]^f \ar@/_1.5pc/[rr]_g &\Downarrow a'&
                Y \ar[r]^w & Z
       }
        }\, ,
\end{align}
then uniqueness of $a'$, together with equation (\ref{antiwhisker_eqn}) gives us
\[
        \raisebox{36pt}{
       \xymatrix{
               &X \ar[dr]^f \\
               W \ar[ur]^{v'} \ar[dr]_{v'} & \Downarrow b & Y \\
               & X \ar[ur]_g
       }
       }
        \equals
        \raisebox{36pt}{
        \xymatrix{
               &&\\
               W\ar[r]^{v'}&X \ar@/^1.5pc/[rr]^f \ar@/_1.5pc/[rr]_g
               &\Downarrow a'
               & Y
       }
        }\, .
\]
This is precisely the diagram (\ref{2cf4.diag}) with $v=\id_X$, $u=v'$, $u'=\id_W$ and 
$\epsilon$ the identity 2-arrow.
%
Hence 2CF4 holds.\endofproof
\end{itemize}

The proof of theorem \ref{bicat_frac_exists} is written using only the language of
2-categories, so can be generalised from $C$ to other 2-categories. This approach
will be taken up in \cite{Roberts2}.

The second main result of the paper is that we want to know when this
bicategory of fractions is equivalent to a bicategory of anafunctors, as the latter 
bicategory has a much simpler construction.

\begin{theorem}\label{anafunctors_localise}
Let $(S,J)$ be a subcanonical unary site and let the full sub-2-category 
$C\into \Cat(S)$ admit base change along arrows in $J$.
Then there is an equivalence of bicategories
\[
       C_\ana(J) \simeq C[W_J^{-1}]
\]
under $C$.
\end{theorem}

\proof
Let us show the conditions in proposition \ref{comparison_thm} hold. To begin with, the 2-functor
$\alpha_J\colon C \to C_{ana}(J)$ sends $J$-equivalences to equivalences by 
proposition \ref{W-inverting_alpha}.
\begin{itemize}
\item[EF1.] $\alpha_J$ is the identity on 0-cells, and hence surjective on objects.

\item[EF2.] This is equivalent to showing that for any anafunctor $(U,f)\colon X\gento Y$
there  are functors $w,g$ such that $w$ is in $W_J$ and
\[
(U,f) \stackrel{\sim}{\Rightarrow} \alpha_J(g)\circ\alpha_J(w)^{-1}
\]
 where $\alpha_J(w)^{-1}$ is some pseudoinverse for $\alpha_J(w)$.

 Let $w$ be the functor $X[U] \to X$ and let $g=f\colon X[U] \to Y$. First, note that
 \[
       \xymatrix{
               & \ar[dl] X[U] \ar[dr]^= &\\
               X       &&              X[U]
       }
 \]
 is a pseudoinverse for
 \[
        \alpha_J(w)
        \equals
        \left(\raisebox{24pt}{
       
       	\xymatrix{
               & \ar[dl]_{=} X[U][U] \ar[dr] &\\
               X[U]     &&                     X
       	}
        
       }\right)\,.
 \]
Then the composition $ \alpha_J(f)\circ\alpha_J(w)^{-1}$ is
\[
       \xymatrix{
               & \ar[dl] X[U\times_U U \times_U U]\ar[dr]\\
               X                               &&                       Y\; ,
       }
\]
which is just $(U,f)$ (recall we have the equality $U\times_U U \times_U U = U$ 
by remark \ref{remark:pullback_of_x_along_x}).

\item[EF3.] If $a\colon(X_0,f)\Rightarrow(X_0,g)$ is a transformation of anafunctors for
functors $f,g\colon X\to Y$, it is given by a natural transformation 
\[
       f \Rightarrow g\colon X = X[X_0 \times_{X_0} X_0]  \to Y.
\]
Hence we get a unique natural
transformation $a\colon f\Rightarrow g$ such that $a$ is the image of $a'$ under
$\alpha_J$.
\endofproof
\end{itemize}

We now give a series of results following from this theorem, using basic properties of
pretopologies from section \ref{sites_categories}. 

\begin{corollary}\label{equivalent_anafunctors}
When $J$ and $K$ are two subcanonical singleton pretopologies on $S$ such that
$J_{un}=K_{un}$, for example $J$ cofinal in $K$, there is an equivalence of bicategories
\[
       C_\ana(J) \simeq C_\ana(K).
\]
\end{corollary}

The class of maps in $\Top$ of the form $\coprod U_i \to X$ for an
open cover $\{U_i\}$ of $X$ form a singleton pretopology. This is because $\mathcal{O}$
is a \emph{superextensive} pretopology (see the appendix).
Given a site with a superextensive pretopology $J$,
we have the following result which is useful when $J$ is not a singleton
pretopology (the \emph{singleton} pretopology $\amalg J$ is defined analogously 
to the case of $\Top$, details are in the appendix).

\begin{corollary}
Let $(S,J)$ be a superextensive site where $J$ is a subcanonical pretopology. Then
\[
       C[W_{J_{un}}^{-1}] \simeq C_\ana(\amalg J).
\]
\end{corollary}
\proof
This essentially follows by lemma \ref{J-coprodJ-epis}.
\endofproof

Obviously this can be combined with previous results, for example if $K$ is cofinal in $\amalg J$,
 for $J$ a non-singleton pretopology, $K$-anafunctors localise $C$ at the
class of $J_{un}$-equivalences.

Finally, given WISC we have a bound on the size of the hom-categories, up to equivalence.

\begin{theorem}
Let $(S,J)$ be a subcanonical unary site satisfying WISC with $S$ locally small
and let $C\into \Cat(S)$ admit base change along arrows in $J$.
Then any localisation $C[W_J^{-1}]$ is locally essentially small. 
\end{theorem}

Recall that this localisation can be chosen such that the class of objects is the same as the class
of objects of $C$, and so it is not necessary to consider additional set-theoretic mechanisms for dealing with
large (2-)categories here.

We note that the issue of size of localisations is not touched on in \cite{Pronk_96}. even 
though such issues are commonly addressed in localisation of 1-categories. If we have 
a specified bound on the hom-sets of $S$ and also know that some WISC${}_\kappa$ holds,
then we can put specific bounds on the size of the hom-categories of the localisation.
This is important if examining fine size requirements or implications for localisation 
theorems such as these, for example higher  versions of locally presentable categories.

\section{Examples}\label{examples}

The simplest example is when we take the trivial singleton pretopology $triv$, where
covering families are just single isomorphisms: $triv$-equivalences are 
internal equivalences and, up to equivalence, localisation at $W_{triv}$ does nothing. It is
worth pointing out that if we localise at $W_{triv_{un}}$, which is equivalent to 
considering anafunctors with source leg having a split epimorphism for its
object component, then by corollary \ref{equivalent_anafunctors} this is equivalent
to localising at $W_{triv}$, so $C_{ana}(triv_{un}) \simeq C_{ana}(triv)\simeq C$.

The first non-trivial case is that of a regular category with 
the canonical singleton pretopology $\can$. This is the setting of \cite{Bunge-Pare_79}.
Recall that $W_J^{BP}$ is the class of Bunge-Par\'e $J$-equivalences (definition
\ref{def:PB_weak_equiv}). For now, let $C$ denote either $\Cat(S)$ or $\Gpd(S)$.

\begin{proposition}
Let $(S,J)$ be a finitely complete unary site with $J$ saturated. Then we have
\[
C[(W_J^{BP})^{-1}] \simeq C[W_J^{-1}]
\]
\end{proposition}

This is merely a restatement of the fact Bunge-Par\'e $J$-equivalences
and ordinary $J$-equivalences coincide in this case.

\begin{corollary}\label{PB_are_the_same_weak_equivs}
The canonical singleton pretopology $\can$ on a finitely complete category $S$
is saturated. Hence $W_{\can}^{BP} = W_{\can}$ for this site, and 
\[
C[(W_{\can}^{BP})^{-1}] \simeq C[W_{\can}^{-1}]\simeq C_\ana(\can)
\]
\end{corollary}

We can combine this corollary with corollary \ref{equivalent_anafunctors} so that the 
localisation of either $\Cat(S)$ or $\Gpd(S)$ at the Bunge-Par\'e weak equivalences 
can be calculated using $J$-anafunctors for $J$ cofinal in $\can$. We note that $\can$
does not satisfy WISC in general (see proposition \ref{U-small_G-sets} and the comments
following), so the localisation might not be locally essentially small.

The previous corollaries deal with the case when we are interested in the 2-categories
consisting of all of the internal categories or groupoids in a site. However, for many applications
of internal categories/groupoids it is not sufficient to take all of $\Cat(S)$ or $\Gpd(S)$.
One widely used example is that of Lie groupoids, which are groupoids internal to the
category of (finite-dimensional) smooth manifolds such that source and target maps 
are submersions (more on these below). Other examples are used in the theory
of algebraic stacks, namely
groupoids internal to schemes or algebraic spaces. Other types of such \emph{presentable}
stacks use groupoids internal to some site with specified conditions on the source and target
maps. Although it is not covered explicitly in the literature, it is possible to consider presentable
stacks of categories, and this will be taken up in future work \cite{Roberts1}.

We thus need to furnish examples of sub-2-categories $C$, specified by restricting
the sort of maps that are allowed for source and target, that admit base change
along some class of arrows. The following lemma gives a sufficiency condition for this
to be so.

\begin{lemma}\label{lemma:existence_of_base_change}
Let $\Cat^\mathcal{M}(S)$ be defined as the full sub-2-category of $\Cat(S)$ with
objects those categories such that the source and target maps belong to a singleton pretopology
$\mathcal{M}$. Then $\Cat^\mathcal{M}(S)$ admits base change along arrows in $\mathcal{M}$, 
as does the corresponding 2-category $\Gpd^\mathcal{M}(S)$ of groupoids.
\end{lemma}
\proof
Let $X$ be an object of $\Cat^\mathcal{M}(S)$ and $f\colon M\to X_0 \in \mathcal{M}$. In the
following diagram, all the squares are pullbacks and all arrows are in $\mathcal{M}$.
\[
\SelectTips {cm}{}%
\xymatrix{
       X[M]_1 \ar[d] \ar[r] \ar @/_2.4pc/  [dd]_{s'} \ar @/^1pc/[rr]^{t'} & X_1\times_{X_0} M \ar[r] \ar[d] & M \ar[d]  \\
       M\times_{X_0} X_1 \ar[d] \ar[r] & X_1 \ar[r] \ar[d] & X_0 \\
       M \ar[r] & X_0 
}
\]
The maps marked $s',t'$ are the source and target maps for the base change along
$f$, so $X[M]$ is in $\Cat^\mathcal{M}(S)$. The same argument holds for groupoids verbatim.
\endofproof

In practice one often only wants base change along a subclass of $\mathcal{M}$,
such as the class of open covers sitting inside the class of open maps in $\Top$. 
We can then apply theoerems \ref{bicat_frac_exists} and \ref{anafunctors_localise}
to the 2-categories $\Cat^\mathcal{M}(S)$ 
and $\Gpd^\mathcal{M}(S)$ with the classes of $\mathcal{M}$-equivalences, and indeed
to sub-2-categories of these, as we shall in the examples below.

We shall focus of a few concrete cases to show how the results of this paper subsume
similar results in the literature proved for specific sites.

The category of smooth manifolds is not finitely complete so the localisation results in this 
section so far do not apply to it. There are two ways around this. The first is to expand the category
of manifolds to a category of smooth spaces which \emph{is} finitely complete 
(or even cartesian closed). In that
case all the results one has for finitely complete sites can be applied. The other
is to take careful note of which finite limits are actually needed, and show that all constructions
work in the original category of manifolds. There is then a hybrid approach, which is to
work in the expanded category, but point out which results/constructions actually fall inside
the original category of manifolds. Here we shall take the second approach. First, let us pin down
some definitions.
\begin{definition}
Let $\Diff$ be the category of smooth, finite-dimensional manifolds. A \emph{Lie category} is
a category internal to $\Diff$ where the source and target maps are submersions (and hence 
the required pullbacks exist). A \emph{Lie groupoid} is a Lie category which is a groupoid.
A \emph{proper} Lie groupoid is one where the map $(s,t)\colon X_1 \to X_0 \times X_0$ is
proper. An \emph{\'etale} Lie groupoid is one where the source and target maps are local
diffeomorphisms.
\end{definition}

By lemma \ref{lemma:existence_of_base_change} the 2-categories of Lie categories, 
Lie groupoids and proper Lie groupoids admit base change along any of the following 
classes of maps: open covers ($\amalg\mathcal{O}$), surjective local diffeomorphisms ($\acute{e}t$), 
surjective submersions ($Subm$).
The 2-categories of \'etale Lie groupoids and proper \'etale Lie groupoids admit base change
along arrows in $\acute{e}t$ and $Subm$. We should note that
we have $\amalg\mathcal{O}$ cofinal in $\acute{e}t$, which is cofinal in
$Subm$.

We can thus apply the main results of this paper to the sites $(\Diff,\mathcal{O})$, 
$(\Diff,\amalg\mathcal{O})$, $(\Diff,\acute{e}t)$
and $(\Diff,Subm)$ and the 2-categories of Lie categories, Lie groupoids, proper Lie goupoids
and so on. However, the definition of weak equivalence we have here, involving $J$-locally 
split functors, is not one that apppears in the Lie groupoid literature, which is actually Bunge-Par\'e
$Subm$-equivalence. However, we have the following result:
\begin{proposition}\label{Subm-equiv_are_BP-equiv}
A functor $f\colon X\to Y$ between Lie categories is a $Subm$-equivalence if and only if 
it is a Bunge-Par\'e
$Subm$-equivalence.
\end{proposition}
Before we prove this, we need a lemma proved by Ehresmann.
\begin{lemma}{\cite{Ehresmann_59}}
For any Lie category $X$, the subset of invertible arrows, $X_1^{iso} \into X_1$ is an open
submanifold. 
\end{lemma}
Hence there is a Lie groupoid $X^{iso}$ and an identity-on-objects functor $X^{iso} \to X$
which is universal for functors from Lie groupoids. In particular, a natural isomorphism between
functors with codomain $X$ is given by a component map that factors through $X_1^{iso}$, and 
the induced source and target maps $X_1^{iso} \to X_0$ are submersions.

\proof (proposition \ref{Subm-equiv_are_BP-equiv})
Full faithfulness is the same for both definitions, so we just need to show that $f$ is $Subm$-locally
split if and only if it is essentially $Subm$-surjective. We first show the forward implication.

The special case of a $\amalg\mathcal{O}$-equivalence between Lie groupoids is a small
generalisation of the proof of proposition 5.5 in \cite{Moerdijk-Mrcun_03}, which states than 
an internal equivalence of Lie groupoids is a Bunge-Par\'e $Subm$-equivalence. Since $\amalg\mathcal{O}$ is 
cofinal in $Subm$, a $Subm$-equivalence is a $\amalg\mathcal{O}$-equivalence, hence a
Bunge-Par\'e $Subm$-equivalence.

For the case when $X$ and $Y$ are Lie categories, we use the fact that we can define 
$X_0\times_{Y_0}Y_1^{iso}$ and that the local sections constructed in Moerdijk-Mr\v cun's proof
factor through this manifold to set up the proof as in the groupoid case.

For the reverse implication, the construction in the first half of the proof of proposition 
\ref{BP_equiv_iff_weak_equiv}  goes through verbatim, as all the pullbacks used 
involve submersions.
\endofproof

The need to localise the category of Lie groupoids at $W_{Subm}$ 
was perhaps first noted in \cite{Pradines_89},
where it was noted that something other than the standard construction of a category of fractions 
was needed. However Pradines lacked the necessary 2-categorical localisation results.
Pronk considered the sub-2-category of \'etale Lie groupoids, also localised at $W_{Subm}$, in order
to relate these groupoids to differentiable \'etendues \cite{Pronk_96}. Lerman discusses the 2-category
of orbifolds \emph{qua} stacks \cite{Lerman_10} and argues that it should be a localisation 
of the 2-category of proper \'etale
Lie groupoids (again at $W_{Subm}$). These three cases use different constructions of the 2-categorical
localisation: Pradines used what he called \emph{meromorphisms}, which are equivalence classes of 
butterfly-like diagrams and are related to Hilsum-Skandalis morphisms, Pronk introduces the techniques 
outlined in this paper, and Lerman uses Hilsum-Skandalis morphisms, also known as right principal 
bibundles. 

Interestingly, \cite{Colman_10} considers this localisation of the 2-category of Lie groupoids then
considers a further localisation, not given by the results of this paper.\footnote{In fact this is the only
2-categorical localisation result involving internal categories or groupoids known to the author to 
\emph{not} be covered by theorem \ref{bicat_frac_exists} or its sequel \cite{Roberts2}.} 
Colman in essence shows that the full
sub-2-category of topologically discrete groupoids, i.e. ordinary small groupoids, is a localisation at 
those internal functors which induce an equivalence on fundamental groupoids.

Our next example is that of topological groupoids, which correspond to various flavours of stacks on 
the category $\Top$. The idea of weak equivalences of topological groupoids predates the 
case of Lie groupoids, and \cite{Pradines_89} credits it to Haefliger, van Est,
and \cite{Hilsum-Skandalis_87}. In particular the first two were ultimately interested in
defining the fundamental group of a foliation, that is to say, of the topological groupoid
associated to a foliation, considered up to weak eqivalence. 

However more recent examples have focussed on topological stacks, or variants thereon.
In particular, in parallel with the algebraic and differentiable cases, the topological stacks
for which there is a good theory correspond to those topological groupoids with conditions
on their source and target maps. Aside from \'etale topological groupoids (which were considered
by \cite{Pronk_96} in relation to \'etendues), the real advances here have come from
work of Noohi, starting with \cite{Noohi_05a}, who axiomatised the 
concept of \emph{local fibration} and asked that the source and target maps of topological
groupoids are local fibrations.

\begin{definition}
A singleton pretopology $LF$ in $\Top$ is called a class of \emph{local fibrations} if the following 
conditions hold:\footnote{We have packaged the conditions in a way slightly different to 
\cite{Noohi_05a}, but the definition is in fact identical.}
\begin{enumerate}
\item $LF$ contains the open embeddings
\item $LF$ is stable under coproducts, in the sense that $\coprod_{i\in I} X_i \to Y$ is in $LF$
if each $X_i\to Y$ is in $LF$
\item $LF$ is local on the target for the open cover pretopology. 
That is, if the pullback of a map $f\colon X\to Y$ along an open cover of $Y$ is in $LF$, then $f$
is in $LF$.
\end{enumerate}
\end{definition}
Conditions 1. and 2. tell us that $\amalg\mathcal{O} \subset LF$, and that $LF$ is $\amalg J$
for some superextensive pretopology $J$ containing the open embeddings as singleton 
`covering' families (beware the misleading terminology here: covering families are not
assumed to be jointly surjective). Note that $LF$ will not be subcanonical,
by condition 1. As an example, given any of the following pretopologies $K$:
\begin{itemize}
\item Serre fibrations,
\item Hurewicz fibrations,
\item open maps,
\item split maps,
\item projections out of a cartesian product,
\item isomorphisms;
\end{itemize}
one can define a class of local fibrations by choosing those maps which are in $K$ on pulling
back to an open cover of the codomain. Such maps are then called \emph{local $K$}. 
As an example of the usefulness of this concept, the topological stacks corresponding to
topological groupoids with local Hurewicz fibrations as source and target have a nicely behaved
homotopy theory. The case of \'etale groupoids corresponds to the last named class of maps,
which give us local isomorphisms, i.e.~\'etale maps.
We can then apply lemma \ref{lemma:existence_of_base_change} and theorem \ref{bicat_frac_exists}
to the 2-category $\Grp^{LF}(\Top)$ to localise at the class $W_{\amalg \mathcal{O}}$ (as 
$\amalg \mathcal{O} \subset LF$), or any other singleton pretopology contained in $LF$,
using anafunctors whenever this pretopology is subcanonical. Note that if $C$ satisfies WISC,
so will the corresponding $LF$, although this is probably not necessary to consider in the presence
of full AC.

A slightly different approach is taken in \cite{Carchedi_12}, where the author introduces a new
pretopology on the category $CGH$ of compactly generated Hausdorff spaces. We give a definition
equivalent to the one in \emph{loc cit}.
\begin{definition}
A (not necessarily open) cover $\{V_i\into X\}_{i\in I}$ is called a $\mathcal{CG}$-cover if for any map
$K\to X$ from a compact space $K$, there is a finite open cover $\{U_j \into K\}$ which refines 
the cover $\{V_i\times_X K\to K\}_{i\in I}$. $\mathcal{CG}$-covers form a pretopology $\mathcal{CG}$
on $CGH$.
\end{definition}
Compactly generated stacks then correspond to groupoids in $CGH$ such that source and
target maps are in the pretopology $\mathcal{CG}_{un}$. Again, we can localise 
$\Gpd^{\mathcal{CG}}(CGH)$ at $W_{\mathcal{CG}_{un}}$ using lemma 
\ref{lemma:existence_of_base_change} and theorem \ref{bicat_frac_exists}, and anafunctors
can be again pressed into service.

We now arrive at the more involved case of algebraic stacks (cf.~the 
continually growing \cite{StacksProj}
for the extent of the theory of algebraic stacks), which were the first presentable
stacks to be defined. There are some subtleties about the site of definition for algebraic stacks,
and powerful representability theorems, but we can restrict to three main cases: groupoids
in the category of affine schemes $\Aff = \Ring^{op}$; groupoids in the category $\Sch$ of schemes; and
groupoids in the category $\AlgSp$ of algebraic spaces. Algebraic spaces reduce to
algebraic stacks on $\Sch$ represented by groupoids with trivial automorphism groups, and the category of
schemes is a subcategory of $Sh(\Aff)$, so we shall just consider the case when our ambient category
is $\Aff$. In any case, all the special properties of classes of maps in all three sites are ultimately
defined in terms of properties of ring homomorphisms.
Note that groupoids in $\Aff$ are exactly the same thing as cogroupoid objects in $\Ring$, which are more 
commonly known as \emph{Hopf algebroids}. 

Despite the possibly unfamiliar language used by algebraic geometry, algebraic stacks reduce 
to the following semiformal definition. We fix three singleton pretopologies on our site $\Aff$: $J$, $E$ 
and $D$ such that $E$ and $D$ are local on the target for the pretopology $J$. An algebraic 
stack then is a stack on $\Aff$ for the pretopology $J$ which `corresponds' to a groupoid $X$ in $\Aff$
 such that source and target maps belong to $E$
and $(s,t)\colon X_1 \to X_0^2$ belongs to $D$. We recover the algebraic stacks by localising the 2-category
of such groupoids at $W_E$ (this claim of course needs substantiating, something we will not do here
for reasons of space, referring rather to \cite{Pronk_96,Schappi_12} and the forthcoming \cite{Roberts1}).

In practice, $D$ can be something like closed maps (to recover Hausdorff-like conditions) or all maps,
and $E$ consists of either smooth or \'etale maps, corresponding to Artin and Deligne-Mumford stacks
respectively. $J$ is then something like the \'etale topology (or rather, the singleton pretopology associated
to it, as the \'etale topology is superextensive), and we can apply lemma \ref{lemma:existence_of_base_change}
to see that base change exists along $J$, along with the fact that asking for $(s,t) \in D$ is automatically
stable under forming the base change. In practice, a variety of combinations of $J,E$ and $D$ are used,
as well as passing from $\Aff$ to $\Sch$ and $\AlgSp$, so there are various compatibilities to check in order
to know one can apply theorem \ref{bicat_frac_exists}.

A final application we shall consider is when our ambient category consists of algebraic objects. As 
mentioned in section 2, a number of authors have considered localising groupoids in 
Mal'tsev, or Barr-exact, or protomodular, or semi-abelian categories, which are hallmarks of categories
of algebraic objects rather than spatial ones, as we have been considering so far.

In the case of groupoids in $\Grp$ (which, as in any Mal'tsev category, coincide with the internal categories)
it is a well-known result that they can be described using \emph{crossed modules}. 
\begin{definition}
A \emph{crossed module} (in $\Grp$) is a homomorphism $t\colon G\to H$ together with a homomorphism
$\alpha\colon H\to \Aut(G)$ such that $t$ is $H$-equivariant (using the conjugation action of $H$ on itself),
and such that the composition $\alpha\circ t\colon G\to\Aut(G)$ is the action of $G$ on itself by conjugation.
A crossed module is often denoted, when no confusion will arise, by $(G\to H)$.
A morphism $(G \to H) \to (K\to L)$ of crossed modules is a pair of maps $G\to K$ and $H\to L$
making the obvious square commute, and commuting with all the action maps.
\end{definition}

Similar definitions hold for groups internal to cartesian closed categories, and even just finite-product
categories if one replaces $H\to \Aut(G)$ with its transpose $H\times G\to G$. Ultimately of course there
is a definition for crossed modules in semiabelian categories (e.g.~\cite{AMMV_10}), but we shall consider just groups.
There is a natural definition of 2-arrow between maps of crossed modules, but the specifics are not
important for the present purposes, so we refer to \cite[definition 8.5]{Noohi_05c} for details. 
The 2-categories of groupoids internal to $\Grp$ and crossed modules are equivalent, so we shall just 
work with the terminology of the latter.

Given the result that crossed modules correspond to pointed, connected homotopy 2-types,
it is natural to ask if all maps of such arise from maps between crossed modules. The answer is, 
perhaps unsurprisingly, no, as one needs maps which only \emph{weakly} preserve the group structure. 
One can either write down the definition of some generalised form of map (\cite[definition 8.4]{Noohi_05c}),
or localise the 2-category of crossed modules (\cite{Noohi_05c} considers a model structure on the 
category of crossed modules). To localise the 2-category of crossed modules we can consider
the singleton pretopology $epi$ on $\Grp$ consisting of the epimorphisms, and localise $\Gpd(\Grp)$ at $W_{epi}$.

There are potentially interesting sub-2-categories of crossed modules that one might want to consider,
for example, the one corresponding to \emph{nilpotent} pointed connected 2-types. These are crossed
modules $t\colon G \to H$ where the cokernel of $t$ is a nilpotent group and the (canonical) action of $\coker t$
on $\ker t$ is nilpotent. The correspondence between such crossed modules and the corresponding internal
groupoids is a nice exercise, as well as seeing that this 2-category admits base change for the pretopology $epi$.

\appendix
\section{Superextensive sites}

The usual sites of topological spaces, manifolds and schemes all share a common property:
one can (generally) take coproducts of covering families and end up with a cover.
In this appendix we gather some results that generalise this fact, none of which are especially 
deep, but help provide examples of bicategories of anafunctors. Another reference for superextensive
sites is \cite{Shulman_12}.

\begin{definition}{\cite{Carboni-Lack-Walters_93}}
A \emph{finitary} (resp. \emph{infinitary}) \emph{extensive} category is a category
with finite (resp. small) coproducts such that the following condition holds: let $I$ be a
a finite set (resp. any set), then, given a collection of commuting diagrams
\[
       \xymatrix{
               X_i \ar[r] \ar[d] &Z \ar[d] \\
               A_i \ar[r] & \coprod_{i\in I} A_i\;,
       }
\]
one for each $i\in I$, the squares are all pullbacks if and only if the collection $\{X_i
\to Z\}_{i\in I}$ forms a coproduct diagram.
\end{definition}

In such a category there is a strict initial object: given a map $A \to 0$ with $0$
initial, we have $A \simeq 0$.

\begin{example}
       $\Top$ is infinitary extensive. $\Ring^{op}$, the category of affine schemes,
       is finitary extensive.
\end{example}

\begin{example}

\end{example}

In $\Top$ we can take an open cover $\{U_i\}_I$ of a space $X$ and replace it with
the single map $\coprod_I U_i \to X$, and work just as before using this new sort of
cover, using the fact $\Top$ is extensive. The sort of sites that mimic this behaviour
are called \emph{superextensive}.

\begin{definition}{(Bartels-Shulman)}
A \emph{superextensive site} is an extensive category $S$ equipped with a
pretopology $J$ containing the families
\[
       (U_i \to \coprod_I U_i)_{i\in I}
\]
and such that all covering families are bounded; this means that for a finitely
extensive site, the families are finite, and for an infinitary site, the families are small.
The pretopology in this instance will also be called superextensive.
\end{definition}

\begin{example}
Given an extensive category $S$, the \emph{extensive pretopology} has as covering
families the bounded collections $(U_i \to \coprod_I U_i)_{i\in I}$. The pretopology on
any superextensive site contains the extensive pretopology.
\end{example}

\begin{example}
The category $\Top$ with its usual pretopology of open covers is a superextensive
site.
\end{example}

\begin{example}
An elementary topos with the coherent pretopology is finitary superextensive, and a Grothendieck
topos with the canonical pretopology is infinitary superextensive.
\end{example}

Given a superextensive site $(S,J)$, one can form the class $\amalg J$ of arrows of the
form $\coprod_I U_i \to A$ for covering families $\{U_i \to A\}_{i\in I}$ in $J$ (more precisely, 
all arrows isomorphic in $S/A$ to such arrows).

\begin{proposition}
The class $\amalg J$ is a singleton pretopology, and is subcanonical if and only if $J$
 is.
\end{proposition}

\proof
Since isomorphisms are covers for $J$ they are covers for $\amalg J$. The pullback
of a $\amalg J$-cover $\coprod_I U_i \to A$ along $B \to A$ is a $\amalg J$-cover as
coproducts and pullbacks commute by definition of an extensive category. Now for
the third condition we use the fact that in an extensive category a map
\[
       f\colon B \to \coprod_I A_i
\]
implies that $B\simeq \coprod_I B_i$ and $f=\coprod_i f_i$. Given $\amalg J$-covers
$\coprod_I U_i \to A$ and $\coprod_J V_j \to (\coprod_I U_i)$, we see that
$\coprod_J V_j \simeq \coprod_I W_i$ for some objects $W_i$. By the previous point, the pullback
\[
       \coprod_I U_k \times_{\coprod_I U_{i'}} W_i
\]
is a $\amalg J$-cover of $U_i$, and hence
$(U_k \times_{\coprod_I U_{i'}} W_i \to U_k)_{i\in I}$
 is a $J$-covering family for each $k\in I$. Thus
\[
       (U_k \times_{\coprod_I U_{i'}} W_i \to A)_{i,k\in I}
\]
is a $J$-covering family, and so
\[
\coprod_J V_j \simeq \coprod_{k\in I} \left( \coprod_{i\in I} U_k \times_{\coprod_I U_{i'}}
W_i\right) \to A
\]
is a $\amalg J$-cover.\\
The map $\coprod_I U_i \to A$ is the coequaliser of $\coprod_{I\times I} U_i \times_A
U_j \st \coprod_I U_i$ if and only if $A$ is the colimit of the diagram in definition
\ref{defn:effective_cov_fam}. Hence $(\coprod_I U_i \to A)$ is effective if and only if $
(U_i \to A)_{i\in I}$ is effective
\endofproof

Notice that the original superextensive pretopology $J$ is generated by the union of
$\amalg J$ and the extensive pretopology.

One reason we are interested in superextensive sites is the following.

\begin{lemma}\label{J-coprodJ-epis}
In a superextensive site $(S,J)$, we have  $J_{un} = (\amalg J)_{un}$.
\end{lemma}

This means we can replace the singleton pretopology $J_{un}$ (e.g. local-section-admitting
maps of topological spaces) with the singleton pretopology $\amalg J$ (e.g. disjoint unions of open covers)
when defining anafunctors. This makes for much smaller pretopologies in practice.

One class of extensive categories which are of particular interest is those that also
have finite/small limits. These are called \emph{lextensive}. For example, $\Top$ is
infinitary lextensive, as is a Grothendieck topos. In contrast, an elementary topos is 
in general only finitary lextensive. We end with a lemma about WISC.

\begin{lemma}
If $(S,J)$ is a superextensive site, $(S,J)$ satisfies WISC if and only if $(S,\amalg J)$
does.
\end{lemma}

One reason for why superextensive sites are so useful is the following result from 
\cite{Schappi_12}.

\begin{proposition}[\cite{Schappi_12}]
Let $(S,J)$ be a superextensive site, and $F$ a stack for the extensive topology on $S$.
Then the associated stack $\widetilde{F}$ on the site $(S,\amalg J)$ is also the associated
stack for the site $(S,J)$.
\end{proposition}

As a corollary, since every weak 2-functor $F\colon S\to \Gpd$ for extensive $S$ 
represented by an internal groupoid is automatically a stack for the extensive topology,
we see that we only need to stackify $F$ with respect to a singleton pretopology on $S$.
This will be applied in \cite{Roberts1}.

\refs

\bibitem[AMMV 2010]{AMMV_10}
O.~Abbad, S.~Mantovani, G.~Metere, and E.M.~Vitale, \emph{Butterflies are
 fractions of weak equivalences}, preprint (2010). Available from 
 \url{http://users.mat.unimi.it/users/metere/}.

\bibitem[Aczel 1978]{Aczel}
P.~Aczel, \emph{The type theoretic interpretation of constructive set
  theory}, Logic Colloquium '77, Stud.~Logic Foundations Math., vol.~96,
  North-Holland, 1978, pp 55--66.
  
\bibitem[Aldrovandi-Noohi 2009]{Aldrovandi-Noohi_09} 
E.~Aldrovandi, B.~Noohi, \emph{Butterflies I: Morphisms of 2-group stacks}, 
Adv.~Math., \textbf{221}, issue 3 (2009), pp 687--773, [arXiv:0808.3627].

\bibitem[Aldrovandi-Noohi 2010]{Aldrovandi-Noohi_10} 
E.~Aldrovandi, B.~Noohi, \emph{Butterflies II: Torsors for 2-group stacks}, 
Adv.~Math., \textbf{225}, issue 2 (2010), pp 922--976, [arXiv:0909.3350].


\bibitem[Bartels 2006]{Bartels}
T.~Bartels, \emph{{Higher gauge theory I: 2-Bundles}}, Ph.D. thesis, University
 of California Riverside, 2006, [arXiv:math.CT/0410328].

\bibitem[Baez-Lauda 2004]{HDA5}
J.~Baez and A.~Lauda, \emph{Higher dimensional algebra {V}: 2-groups}, Theory
 and Application of Categories \textbf{12} (2004), no.~14, pp 423--491.
 
 \bibitem[Baez-Makkai 1997]{Baez-Makkai_emails}
 J.~Baez and M.~Makkai, \emph{Correspondence on the category theory mailing
 list, January 1997}, available at \url{http://www.mta.ca/~cat-dist/catlist/1999/anafunctors}.

\bibitem[B{\'e}nabou 1973]{Benabou_73}
J.~B{\'e}nabou, \emph{Les distributeurs}, rapport 33,
 Universit{\'e} Catholique de Louvain, Institut de Math{\'e}matique Pure et
 Appliqu{\'e}e, 1973.

\bibitem[B{\'e}nabou 1975]{Benabou_75a}
J.~B{\'e}nabou, \emph{Th\'eories relatives \`a un corpus}, C.~R.~Acad.~Sci.~Paris
 S\'er.~A-B \textbf{281} (1975), no.~20, Ai, A831--A834.
 
\bibitem[B{\'e}nabou 1989]{Benabou_89}
J.~B{\'e}nabou, \emph{Some remarks on 2-categorical algebra}, Bulletin de la
Soci{\'e}t{\'e} Math{\'e}matique de Belgique \textbf{41} (1989), pp 127--194.

\bibitem[B{\'e}nabou 2011]{Benabou_email}
J.~B{\'e}nabou, \emph{Anafunctors versus distributors}, email to Michael Shulman, 
posted on the category theory mailing list 22 January 2011, available from
\url{http://article.gmane.org/gmane.science.mathematics.categories/6485}.

\bibitem[van den Berg 2012]{vdBerg_12}
B.~van~den Berg, \emph{Predicative toposes}, preprint (2012), [arXiv:1207.0959].


\bibitem[Breckes 2009]{Breckes_09}
M.~Breckes, \emph{Abelian metamorphosis of anafunctors in butterflies},
preprint (2009).

\bibitem[Bunge-Par\'e 1979]{Bunge-Pare_79}
M.~Bunge and R.~Par{\'e}, \emph{Stacks and equivalence of indexed categories},
Cah.~Topol. G\'eom.~Diff\'er.~\textbf{20} (1979), no.~4,
pp 373--399.

\bibitem[CLW 1993]{Carboni-Lack-Walters_93}
A.~Carboni, S.~Lack, and R.~F.~C. Walters, \emph{Introduction to
  extensive and distributive categories}, J.~Pure Appl.~Algebra \textbf{84}
  (1993), pp 145--158.

\bibitem[Carchedi 2012]{Carchedi_12}
D.~Carchedi, \emph{Compactly generated stacks: a cartesian-closed theory of
 topological stacks}, Adv. Math. \textbf{229} (2012), no.~6, pp 3339--33397, [arXiv:0907.3925]

\bibitem[Colman 2010]{Colman_10}
H.~Colman, \emph{On the homotopy 1-type of {L}ie groupoids}, Appl.~Categ.~ 
Structures \textbf{19}, Issue 1 (2010), pp 393--423, [arXiv:math/0612257].

\bibitem[Colman-Costoya 2009]{Colman-Costoya_09}
H.~Colman and C.~Costoya, \emph{A {Q}uillen model structure for orbifolds},
  preprint (2009). Available from \url{http://faculty.ccc.edu/hcolman/}.

%

\bibitem[Dwyer-Kan 1980a]{Dwyer-Kan_80a}
W.~G.~Dwyer and D.~M.~Kan, \emph{Simplicial localizations of categories}, J.~Pure
Appl.~Algebra \textbf{17} (1980), no.~3, pp 267--284.
%
%

\bibitem[Dold 1963]{Dold_63}
A.~Dold, \emph{Partitions of unity in the theory of fibrations}, Ann.~Math.~
 \textbf{78} (1963), no.~2, pp 223--255.

\bibitem[Ehresmann 1959]{Ehresmann_59}
C.~Ehresmann, \emph{Cat\'egories topologiques et cat\'egories
 diff\'erentiables}, Colloque G\'eom. Diff.~Globale (Bruxelles, 1958), Centre
 Belge Rech. Math., Louvain, 1959, pp 137--150.

\bibitem[Ehresmann 1963]{Ehresmann_63}
C.~Ehresmann, \emph{Cat\'egories structur\'ees}, Annales de l'Ecole Normale et
 Superieure \textbf{80} (1963), pp 349--426.

\bibitem[EKvdL 2005]{Everaert_et_al_05}
T.~Everaert, R.W. Kieboom, and T.~van~der Linden, \emph{Model structures for
 homotopy of internal categories}, Theory and Application of Categories \textbf{15} (2005), no.~3, pp 66--94.

\bibitem[Gabriel-Zisman 1967]{Gabriel-Zisman}
P.~Gabriel and M.~Zisman, \emph{Calculus of fractions and homotopy theory},
 Springer-Verlag, 1967.

\bibitem[Hilsum-Skandalis 1987]{Hilsum-Skandalis_87}
M.~Hilsum and G.~Skandalis, \emph{Morphismes-orient{\'e}s d`epsaces de feuilles
 et fonctorialit{\`e} en th{\'e}orie de {K}asparov (d'apr{\`e}s une conjecture
 d'{A}. {C}onnes)}, Ann.~Sci.~{\'E}cole Norm. Sup. \textbf{20} (1987),
pp 325--390.

\bibitem[Johnstone 2002]{Elephant}
P.~Johnstone, \emph{Sketches of an elephant, a topos theory compendium}, Oxford
 Logic Guides, vol.~43 and 44, The Clarendon Press Oxford University Press,
 2002.

\bibitem[Joyal-Tierney 1991]{Joyal-Tierney_91}
A.~Joyal and M.~Tierney, \emph{Strong stacks and classifying spaces},
 Category theory ({C}omo, 1990), Lecture Notes in Math., vol.~1488, Springer, 1991,
pp 213--236.

\bibitem[Karaglia 2012]{Karaglia_12}
A.~Karaglia, \emph{Embedding posets into cardinals with $DC_{\kappa}$}, preprint
(2012), [arXiv:1212.4396].

\bibitem[Kelly 1964]{Kelly_64}
G.~M.~Kelly, \emph{Complete functors in homology. I. Chain maps and endomorphisms}
Proc.~Cambridge Philos.~Soc.~\textbf{60} (1964), pp 721--735. 



\bibitem[Lerman 2010]{Lerman_10}
E.~Lerman, \emph{Orbifolds as stacks?}, L'Enseign Math. (2) \textbf{56} (2010),
  no.~3-4, pp 315--363, [arXiv:0806.4160].
  
\bibitem[Lurie 2009a]{Lurie_HTT}
J.~Lurie, \emph{Higher Topos Theory}, Annals of Mathematics Studies \textbf{170},
Princeton University Press, 2009. Available from \url{http://www.math.harvard.edu/~lurie/}.


\bibitem[Mac Lane-Moerdijk 1992]{MacLane-Moerdijk}
S.~MacLane and I.~Moerdijk, \emph{Sheaves in geometry and logic},
 Springer-Verlag, 1992.

\bibitem[Makkai 1996]{Makkai}
M.~Makkai, \emph{Avoiding the axiom of choice in general category theory}, J.~Pure Appl.~Algebra
 \textbf{108} (1996), pp 109--173. Available from \url{http://www.math.mcgill.ca/makkai/}.
 
 \bibitem[MMV 2012]{MMV2012}
 S.~Mantovani, G.~Metere and E.~M.~Vitale, 
 \emph{Profunctors in Mal'cev categories and fractions of functors}, preprint (2012).
 Available from \url{http://perso.uclouvain.be/enrico.vitale/research.html}.




\bibitem[Moerdijk-Mr\v cun 2003]{Moerdijk-Mrcun_03}
I.~Moerdijk and J.~Mr{\v c}un, \emph{Introduction to foliations and lie
 groupoids}, Cambridge studies in advanced mathematics, vol.~91, Cambridge
 University Press, 2003.

 
\bibitem[Mr{\v{c}}un 2001]{Mrcun_01}
 J.~Mr\v{c}un, \emph{The Hopf algebroids of functions on \'etale groupoids and their
 principal Morita equivalence},
J.~Pure Appl.~Algebra \textbf{160} (2001), no.~2-3, pp 249--262. 

\bibitem[Noohi 2005a]{Noohi_05a}
B.~Noohi, \emph{Foundations of topological stacks {I}}, preprint (2005),
 [arXiv:math.AG/0503247]. 

\bibitem[Noohi 2005b]{Noohi_05b}
B.~Noohi, \emph{On weak maps between 2-groups}, preprint (2005), [arXiv:math/0506313].

\bibitem[Noohi 2005c]{Noohi_05c}
B.~Noohi, \emph{Notes on 2-groupoids, 2-groups and crossed-modules}, preprint (2005) 
[arXiv:math/0512106].

\bibitem[Pradines 1989]{Pradines_89}
J.~Pradines, \emph{Morphisms between spaces of leaves viewed as fractions},
Cah.~Topol. G\'eom.~Diff\'er.~Cat\'eg.~\textbf{30} (1989),
 no.~3, pp 229--246, [arXiv:0803.4209].

\bibitem[Pronk 1996]{Pronk_96}
D.~Pronk, \emph{Etendues and stacks as bicategories of fractions}, Compositio
 Math.~\textbf{102} (1996), no.~3, pp 243--303.

\bibitem[Roberts 2013]{Roberts_13}
D.~M.~Roberts, \emph{Con(ZF+ $\neg$WISC)}, preprint, (2013).

\bibitem[Roberts A]{Roberts1}
D.~M.~Roberts, \emph{All presentable stacks are stacks of anafunctors}, forthcoming (A).

\bibitem[Roberts B]{Roberts2}
D.~M.~Roberts, \emph{Strict 2-sites, $J$-spans and localisations}, forthcoming (B).

\bibitem[Sch{\"a}ppi 2012]{Schappi_12}
D.~Sch\"appi, \emph{A characterization of categories of coherent sheaves of certain algebraic stacks}, 
preprint (2012), [arXiv:1206.2764].

\bibitem[Shulman 2012]{Shulman_12}
M.~Shulman, \emph{Exact completions and small sheaves}, Theory and Application of 
Categories, \textbf{27} (2012), no.~7, pp 97--173.

\bibitem[Stacks project]{StacksProj}
The Stacks project authors, \emph{Stacks project}, \url{http://stacks.math.columbia.edu}.

\bibitem[Street 1980]{Street_80}
R.~Street, \emph{Fibrations in bicategories},
Cah.~Topol.~G\'eom.~Diff\'er.~Cat\'eg.~\textbf{21} (1980),
pp 111--160.
 
\bibitem[TXL-G 2004]{Tu-Xu-LaurentGengoux_04}
J.-L. Tu, P.~Xu and C.~Laurent-Gengoux
\emph{Twisted K-theory of differentiable stacks},
Ann.~Sci.~\'Ecole Norm.~Sup.~(4) \textbf{37} (2004), no.~6, pp 841--910, [arXiv:math/0306138].

\bibitem[Vitale 2010]{Vitale_10}
E.~M. Vitale, \emph{Bipullbacks and calculus of fractions}, Cah.~Topol. 
G\'eom.~Diff\'er.~Cat\'eg.~\textbf{51} (2010), no.~2, pp 83--113. 
Available from \url{http://perso.uclouvain.be/enrico.vitale/}.

\endrefs

\end{document}